\documentclass[11pt]{article}

\RequirePackage[OT1]{fontenc} \RequirePackage{amsthm,amsmath}

\usepackage{graphics,epsfig,epsf,psfrag}
\usepackage{url}
\usepackage{color, appendix}
\usepackage{graphicx}
\usepackage{amsmath, amsfonts, amsthm,amssymb}
\usepackage{epstopdf}
\usepackage{hyperref}
\usepackage{enumerate, framed}
\usepackage{graphicx,lscape,rotate}
\usepackage{algorithmic, algorithm}
\usepackage{epsfig}
\usepackage{mathrsfs}
\usepackage{amsfonts}
\usepackage{mathtools}
\usepackage{tasks}
\usepackage[nobysame]{amsrefs}
\usepackage{ifthen}
\usepackage{bbm}
\graphicspath{{eps/}}
\usepackage{amsthm}
\usepackage[utf8]{inputenc}
\usepackage[english]{babel}
\usepackage{bm}
\newcommand{\sK}{\mathscr{K}}

\numberwithin{equation}{section}
\theoremstyle{plain}
\newtheorem{theorem}{Theorem}[section]
\newtheorem{assumption}{Assumption}[section]
\newtheorem{remark}{Remark}
\newtheorem{proposition}{Proposition}[section]

\newtheorem{lemma}{Lemma}[section]
\newtheorem{definition}{Definition}
\newtheorem{example}[theorem]{Example}
\theoremstyle{remark}

\newcommand{\F}{ \mathcal{F}}
\newcommand{\R}{ \mathbb{R}}

\newcommand{\s}{\mathcal{L}}

\newcommand{\eq}{\begin{equation}}
\newcommand{\eeq}{\end{equation}}

\begin{document}


\title{Soliton resolution for nonlinear Schr\"odinger type equations in the radial case }
\author{Avy Soffer and Xiaoxu Wu}


\date{}  

\maketitle

\begin{abstract}
\noindent We consider the Schr\"odinger equation with a general interaction term, which is localized in space, for radially symmetric initial data in $n$ dimensions, $n\geq5$. The interaction term may be space-time dependent and nonlinear. Assuming that the solution is bounded in $H^1(\mathbb{R}^n)$ uniformly in time, we prove soliton resolution conjecture (Asymptotic completeness) by demonstrating that the solution resolves into a smooth and localized part and a free radiation in $\s^2_x(\mathbb{R}^n)$ norm. Examples of such equations include inter-critical and super-critical nonlinear Schr\"odinger equations, saturated nonlinearities, time-dependent potentials, and combinations of these terms.
\end{abstract}
{
\hypersetup{linkcolor=black}
\tableofcontents}

\section{Introduction}
\subsection{Background}
A standard approach for studying physical systems described by Hyperbolic and Dispersive PDEs is to use the tools of scattering theory. 
The key starting point is to identify all possible asymptotic states of the system as time approaches positive or negative infinity.
Spectral theory is highly developed for linear systems derived from time-independent Hamiltonians in quantum mechanics. This theory covers both $N$-body scattering and more general linear systems, including linearizations of fluid dynamics, general relativity, and others. Local decay estimates play a central role in this approach, with many references available on this topic, including \cites{AMG1996, FH, ger, KK2014, M1981, SX2016, V2013}. The Enss method, described in \cites{D1980, E1978, P1983}, is also relevant in this context, as are propagation estimates discussed in \cites{HS2000, SS1987, SS1988, Skib} and related works.

However, for time-dependent potentials and nonlinear equations, the scattering problem does not have a similarly general, abstract analysis. Such problems arise in various contexts, including linearizations of nonlinear problems, models for open systems, or as leading approximations to physical systems. See \cite{SW2020} for recent work in this area.

A particular class of problems similar to those treated in this work arises in nonlinear optics, where the governing equations for saturated nonlinearities are used in nonhomogeneous media. This leads to potential terms and coefficients that are both space- and time-dependent. See, for example, \cite{BP3S}.

In this work, which builds on a series of previous works \cites{Sof-W2, Sof-W3, Sof-W5, SW20221} on general scattering theory, we provide definitive answers to a broad class of Nonlinear Schr\"odinger-type equations. We prove that for such equations in five or more dimensions, spherically symmetric initial data with a global solution in $H^1$ leads to an asymptotic state that converges (as expected) to a free wave and a part that is localized. The localized part is concentrated around the origin and is a smooth function of space and time.
We also show, by a simple example how one can recover dispersive estimates, for solutions of the equation via the methods developed in these works. See Section 4.

\subsection{Setting of the problem}
Let $H_0=-\Delta_x$. We consider the general class of Nonlinear Schr\"odinger-type equations of the form:
\begin{equation}
\begin{cases}
    i\partial_t\psi(x,t)=(H_0+\mathcal{N}(|\psi|,|x|,t))\psi(x,t)\\
    \psi(x,0)=\psi_0(x)\in H^1_x(\mathbb{R}^n)\\
    \psi_0\text{ is radial in }x
\end{cases},\quad (x,t)\in \mathbb{R}^{n+1},\quad n\geq 5\label{SE}.
\end{equation}
Here, $\psi(x,t)$ is a complex-valued function of $(x,t)\in \mathbb{R}^{n+1}$, and $\mathcal{N}: \mathbb{R}^3\to \mathbb{R}$, is real-valued. Let $\langle x\rangle:=\sqrt{1+|x|^2}$ and let $\s^p_{x,\sigma}(\mathbb{R}^n)$ denote the weighted $\s^p_x(\mathbb{R}^n)$ space:
\eq
\s^p_{x,\sigma}(\mathbb{R}^n):=\{f(x): \langle x\rangle^{\sigma}f\in \s^p_x(\mathbb{R}^n)\}
\eeq
for $1\leq p\leq \infty$. We impose the uniform decay assumption on $\mathcal{N}(|\psi|, |x|,t)\psi(t)$ and uniform $\s^2_x$ assumption on the interaction $\mathcal{N}$: for some $\sigma>2$,
\eq
\begin{cases}
\mathcal{N}(|\psi|, |x|,t)\psi(t)\in\s^\infty_t\s^q_x(\mathbb{R}^{n+1}),\quad 1 \leq q < 2n/(n+2), \\
 \chi(|x|>1)\mathcal{N}(|\psi|,|x|,t)\psi(t)\in \s^\infty_t\s^{2}_{x,\sigma}(\mathbb{R}^{n+1})\cap\s^\infty_t\s^{q}_x(\mathbb{R}^{n+1})
 \end{cases}.\label{condition}
\eeq
Typical examples are monomial nonlinearities, in particular, some inter-critical and energy super-critical equations of the form:
\eq
\mathcal{N}(|\psi|,|x|,t)=\pm \lambda|\psi|^p\label{texam}
\eeq
for any $p\in (1,\frac{7}{3}), \lambda >0$ when $n=5$, as shown in Example \ref{example}. Here, inter-criticality means that the power of $|\psi|$ in $\mathcal{N}$ is between the mass critical point $\frac{4}{n}$ and energy critical point $\frac{4}{n-2}$. If the interaction contains a supercritical defocusing part, then the energy bound controls also the term
$$
\int |\psi|^{p+2} d^nx.
$$
Using this bound and the assumption of global existence and uniform bound on the $H^1$ norm, we see that the condition\ref{condition} is satisfied for $q=1$ by interpolating the above $L^p$ bound and the $H^1$ bound.\par

When $\mathcal{N}$ is negative, the solution $\psi(t)$ either blows up in a finite time by the Virial identity \cite{Gla}, \cite{O1991}, or exists globally in $\s^2_x(\mathbb{R}^n)$. In this paper, we focus on the asymptotic behavior (as $t\to +\infty$) of solutions that do not blow up in the energy class. More precisely, we assume \emph{a priory}  uniform $H^1_x$ bound:
\eq
\sup\limits_{t\in \mathbb{R}}\|\psi\|_{H_x^1(\mathbb{R}^n)}=E<\infty\label{H1}
\eeq
for some $0<E<\infty$. With the {\em{a priory}} assumption \eqref{H1}, we know that $\psi$ is globally well-posed in $H^1_x(\mathbb{R}^n)$ by iterating the local well-posedness argument, cf. \cite{T2006}. For global existence of general time dependent potentials see e.g., \cite{RSII}. For data whose energy $E$ is small, it is known that the solution $\psi(t)$ eventually scatters to approach a free solution $e^{-itH_0}\psi_+$ in $\s^2_x(\mathbb{R}^n)$ as $t\to +\infty$, see e.g., \cite{Z2006}. For large data, we do not expect this scattering behavior in general, even assuming the energy bound \eqref{H1}. This is due to the existence of soliton solutions. These soliton solutions are generated, for example, by finite-energy solutions $\psi(x)$ to the nonlinear eigenfunction equation:
\eq
-\omega \psi+\Delta_x \psi=F(\psi)\label{sol}
\eeq
for some $\omega>0$. Such solutions are known to be smooth and exponentially rapidly decreasing(see e.g. \cite{C1989}). If one further assumes that $\psi$ is non-negative and spherically symmetric, then there is a unique solution to \eqref{sol} for each $\omega>0$ \cite{C1972}; we refer to this $\psi$ as the {\em{ground state}}. Additionally, there also exist radial solutions that change sign, see \cite{B1983}. We refer to these as \emph{excited states}.

The \emph{Soliton Resolution Conjecture} states that if there is a global solution $\psi(t)$ to the nonlinear Schr\"odinger equation (NLS) \eqref{SE}, then $\psi(t)$ will evolve to a free radiation $e^{-itH_0}\psi_{free}$ plus solitons $\psi_{sol}(t)$, such that
\begin{equation}
\lim\limits_{t\to \infty} \| \psi(t)-e^{-itH_0}\psi_{free}-\psi_{sol}(t)\|_{\s^2_x(\mathbb{R}^n)}=0.
\end{equation}
Knowing all asymptotic states and proving the convergence of the solution to these states is called Asymptotic Completeness (AC).

\subsection{Statement of the result}

In this note, we prove a decomposition of the asymptotic states, with $\psi_{sol}=\psi_{loc}$ first, for the system \eqref{SE}. We assume that \eqref{condition} and \eqref{H1} are valid. Then we show that when $n\geq5$,
\eq
\lim\limits_{t\to \infty}\|\psi(t)-e^{-itH_0}\psi_{free}(x)-\psi_{loc}(t)\|_{\s_x^2(\mathbb{R}^n)}=0\label{goal}
\eeq
holds for some $\psi_{free}\in \s^2_x(\mathbb{R}^n)$ and a localized function $\psi_{loc}(t)$ satisfying
\eq
\|\langle x\rangle^\delta \psi_{loc}(t)\|_{\s^\infty_t\s^2_x(\mathbb{R}^{n+1})}<\infty\label{local}
\eeq
for some $\delta>0$.
\begin{theorem}\label{thm}Let $\psi(t)$ be a solution of the system \eqref{SE} satisfying \eqref{H1}.  If \eqref{condition} is valid for some $\sigma>2$, then \eqref{goal} is true with $\psi_{loc}(t)$ satisfying \eqref{local} when $n\geq 5$.
\end{theorem}
Examples satisfying the conditions of Theorem \ref{thm} include, when $n=5$,
\eq
\mathcal{N}(|\psi|,|x|,t)=\pm \lambda |\psi|^p, \text{ or } W(x,t) \pm \lambda |\psi|^p
\eeq
for all $\lambda>0, p\in (1,4/3],$ and some $W(x,t)$ satisfying that for some $\sigma>2$,
\eq
\begin{cases}
W(x,t)\in \s^\infty_t\s^2_x(\mathbb{R}^{5+1})\\
W(x,t)\in \s^\infty_t\s^2_{\sigma,x}(\mathbb{R}^{5+1})\cap\s^\infty_t\s^\infty_{\sigma,x}(\mathbb{R}^{5+1})
\end{cases},
\eeq
see Examples \ref{example}, \ref{example0}.
\begin{remark}
When we say "inter-critical," we mean that the power of nonlinearity is between the mass critical point ($p=\frac{4}{n}$) and the energy critical point ($p=\frac{4}{n-2}$). Specifically, when $n=5$, the mass critical point is $\frac{4}{n}=\frac{4}{5}$. We note that if $p>\frac{4}{5}$ and is sufficiently close to $\frac{4}{5}$, our method may not work for NLS with $\mathcal{N}=- |\psi(t)|^p$.
\end{remark}
\begin{remark}
Note that Theorem \ref{thm} only requires uniform boundedness when $|x|\geq 1$, which allows for the presence of singularities around the origin at $x=0$. However, if we assume a uniform boundedness condition on the interaction $\mathcal{N}$, we can prove that the localized part is smooth in the sense of \eqref{smooth}.
\end{remark}
Throughout the paper, we use the notation $A\lesssim_s B$ and $A\gtrsim_s B$ to indicate that there exists a constant $C=C(s)>0$ such that $A\leq CB$ and $A\geq CB$, respectively.
When the interaction $\mathcal{N}$ is uniformly bounded in $x$, we can also prove that the localized part is smooth. Specifically, in Theorem \ref{thm2}, we assume that the interaction $\mathcal{N}(|\psi|,|x|,t)\in \s^\infty_t\s^2_x(\mathbb{R}^{n+1})$ satisfies, for all $f\in \s^\infty_tH^1_x(\mathbb{R}^{n+1})$,
\eq
|\mathcal{N}(|f|,|x|,t)|\lesssim_{\|f \|_{\s^\infty_tH^1_x(\mathbb{R}^{n+1})}} \frac{1}{\langle x\rangle^{\sigma}}\quad \text{ for some }\sigma >6:\label{con1}
\eeq
\begin{assumption}\label{asp1} \eqref{con1} holds.
\end{assumption}
Let $V(x,t):=\mathcal{N}(|\psi(t)|,|x|,t)$, $\psi_D(t):=\psi(t)-e^{-itH_0}\psi(0)$, and $V_D(x,t):=\mathcal{N}(|\psi_D(t)|,|x|,t)$. Due to \eqref{H1}, we have $\|\psi_D(t)\|_{\s^\infty_tH^1_x(\mathbb{R}^{n+1})}\leq 2E$.
\begin{assumption}\label{asp4}For all $j=1,2,\cdots,n$, let $P_j:=-i\partial_{x_j}$. One has 
\eq
\| P_j \mathcal{N}(|\psi(t)|,|x|,t)\psi(t)\|_{\s^\infty_t\s^2_{x,2}(\mathbb{R}^{n+1})}\lesssim_E 1.
\eeq
\end{assumption}
\begin{lemma}\label{LempsiD}If \eqref{H1} and \eqref{con1} hold and if Assumption \ref{asp4} is satisfied, then we have
 \eq
 \|\langle x\rangle^{-2}\langle P\rangle^{3/2}\psi_D(t)\|_{\s^\infty_t\s^2_x(\mathbb{R}^{n+1}) }\lesssim_E 1 .
 \eeq
 \end{lemma}
The proof of Lemma \ref{LempsiD} can be found in Section \ref{subsec2.3}.

Based on Lemma \ref{LempsiD}, we know that $\psi(t)\in H^1_x(\mathbb{R}^n)$ and $\langle x\rangle^{-2}\psi_D(t)\in H^{3/2}_x(\mathbb{R}^n)$. Therefore, in the proof of smoothness, we construct $\psi_{loc}(t)$ in terms of $\psi_D(t)$ instead of $\psi(t)$.

We also need the following assumptions related to $\psi_D(t)$:
\begin{assumption}\label{asp2}For some $\sigma>2$ and $\sigma'\in(0,1/2)$ (here we take $\sigma'<1/2$ since we only have limited decay in $|x|$ from $V(x,t)$), we have
\begin{multline}
\|V(x,t)\psi(t)-V_D(x,t)\psi_D(t) \|_{\s^2_{x,\sigma}(\mathbb{R}^n)\cap \s^1_{x,2}(\mathbb{R}^n)}\lesssim_{E} \\\|e^{-itH_0}\psi(0)\|_{\s^2_{x,-\sigma'}(\mathbb{R}^n)}+\|e^{-itH_0}\psi(0)\|_{\s^2_{x,-\sigma'}(\mathbb{R}^n)}^{2/n}\label{con2}.
\end{multline}
\end{assumption}
\begin{remark}Note that $\psi(t)-\psi_D(t)=e^{-itH_0}\psi(0)$ and both $V(x,t)$ and $V_D(x,t)$ are localized functions due to Assumption \eqref{asp1}. That is how we get a factor $\|e^{-itH_0}\psi(0)\|_{\s^2_{x,-\sigma'}(\mathbb{R}^n)}$. In $V(x,t)\psi(t)-V_D(x,t)\psi_D(t)$, there is another term $(V(x,t)-V_D(x,t))\psi_D(t)$. Here we could get a factor $e^{-itH_0}\psi(0) \times \psi_D(t)$ or $\left(e^{-itH_0}\psi(0)\right)^* \times \psi_D(t)$. When $|x|\leq 1$, $\psi_D(t)$ is not uniformly bounded, so in this case, $\s^2_{x,\sigma}(\mathbb{R}^n)$ cannot be controlled by $\|e^{-itH_0}\psi(0)\|_{\s^2_{x,-\sigma'}(\mathbb{R}^n)}$. On the other hand, $$|(V(x,t)-V_D(x,t))\psi_D(t)|\lesssim_E |\psi_D(t)|.$$ Using interpolation, one can obtain, for example,
\eq
\chi(|x|\leq 1)|(V(x,t)-V_D(x,t))\psi_D(t)|\lesssim_E  \chi(|x|\leq 1)|e^{-itH_0}\psi(0) |^{2/n}|\psi_D(t)|.
\eeq
Using that $\sup\limits_{t\in \mathbb{R}}\|\psi_D(t)\|_{\s^{2n/(n-2)}_x(\mathbb{R}^n)}\lesssim \sup\limits_{t\in \mathbb{R}}\|\psi_D(t)\|_{H^{1}_x(\mathbb{R}^n)}\lesssim_E 1$, by H\"older's inequality $(\frac{1}{n}+(\frac{1}{2}-\frac{1}{n})=\frac{1}{2})$, one has
\eq
\begin{split}
\|\chi(|x|\leq 1)(V(x,t)-V_D(x,t))\psi_D(t)\|_{\s^2_{x,\sigma}(\mathbb{R}^n)}\lesssim_E& \| |e^{-itH_0}\psi(0) |^{2/n}\|_{\s^n_{x,-\sigma'}(\mathbb{R}^n)}\|\psi_D(t)\|_{\s^{2n/(n-2)}_x(\mathbb{R}^n)}\\
\lesssim_E& \|e^{-itH_0}\psi(0)\|_{\s^2_{x,-\sigma'}(\mathbb{R}^n)}^{2/n}.
\end{split}
\eeq
\end{remark}
\begin{assumption}\label{asp3}
\eq
\| \langle P\rangle^{3/2}\langle x\rangle^6 V_D(x,t)\psi_D(t)\|_{\s^2_x(\mathbb{R}^n)}\lesssim_{E}1.\label{con5}
\eeq
\end{assumption}
{Under the aforementioned assumptions, if $\mathcal{N}(|\psi(t)|,|x|,t)\in \s^\infty_t\s^2_x(\mathbb{R}^{n+1})$ and \eqref{H1} are satisfied, then \eqref{goal} holds, and
\eq
\| A^2\psi_{loc}(t)\|_{\s^2_{x}(\mathbb{R}^n)}\lesssim_E 1\label{smooth},
\eeq
where $P:=-i\nabla_x$ and $A:=\frac{1}{2}(P\cdot x+x\cdot P)$:
\begin{theorem}\label{thm2}Let $\psi(t)$ denote a solution to the system \eqref{SE} that satisfies \eqref{H1} and \eqref{con1}. If $\mathcal{N}(|\psi|,|x|,t)\in \s^\infty_t\s^2_{x,2}(\mathbb{R}^{n+1})$ satisfies assumptions \ref{asp1}-\ref{asp3}, then \eqref{goal} is true with $\psi_{loc}(t)$ satisfying \eqref{local} and \eqref{smooth}. In particular, when $n\geq 45$, $\sigma> n/2$ in Assumption \ref{asp1}, one has that for some $\delta>1$,
\eq
\| \langle x\rangle^\delta \psi_{loc}(t)\|_{\s^2_x(\mathbb{R}^n)}\lesssim 1 .\label{locallast}
\eeq
\end{theorem}
Typical examples of Theorem \ref{thm2} are
\eq
\mathcal{N}=- \frac{\lambda |\psi|^p}{1+|\psi|^p}\quad \text{ for all }p> 3, \lambda>0, n=5,
\eeq
see example \ref{example2}.

\subsection{ Review of previous works and outline of the proof}
The thrust of this work is the study of general initial data, in particular arbitrary large data.
While it is possible to characterize initial conditions for the solution are understood, even if large, these are all restricted to one channel problems: the solution either scatter to a free wave, or blows-up in a finite time. See e.g. \cites{beceanu2021large,beceanu2018large,beceanupositivity}. Here we focus on the cases which are multi-channel.
Tao's work \cites{T2004, T2007, T2008} showed that the asymptotic decomposition
\eq
\lim\limits_{t\to \infty}\|\psi(t)-e^{-itH_0}\psi_{free}(x)-\psi_{loc}(t)\|_{\dot{H}_x^1(\mathbb{R}^n)}=0
\eeq
holds for NLS with inter-critical nonlinearities in 3 and higher dimensions for radial initial data. When $n\geq 5$, Tao \cite{T2007} extended the decomposition to $H^1_x(\mathbb{R}^n)$ and showed that $\psi_{loc}(t)$ satisfies
\eq
\limsup\limits_{t\to +\infty} \int_{|x|>R} |\psi(x,t)|^2+|\nabla \psi(x,t)|^2 dx\leq c_E(R)
\eeq
for some explicit quantity $c_E(R)$ which goes to zero as $R\to \infty$. To be precise, Tao's method gives a decay rate of $1/\log^cR$ as mentioned in Remark 1.21 in \cite{T2014} (the latest version of \cite{T2007} on the arxiv). For further progress see \cite{killip2021scattering}. Recently, Liu and Soffer \cites{LS2020, LS2021} showed that with radial symmetry assumption, for space dimension $n=3$, if $|\mathcal{N}(|\psi|, |x|,t)|\lesssim \frac{1}{\langle x\rangle^{3+0}} $, then
\eq
\lim\limits_{t\to \infty}\| \psi(t)-e^{-itH_0}\Omega_{+}^*\psi_0-\psi_{wl}(t)\|_{H^1_x(\mathbb{R}^n)}=0\label{eq1}
\eeq
where the weakly localized part $\psi_{wl}$ satisfies
\eq
(\psi_{wl}(t), |x|\psi_{wl}(t))_{\s^2_x(\mathbb{R}^n)}\lesssim |t|^{1/2},\quad t\geq 1
\eeq
and
\eq
\| (x\cdot \nabla_x)^k\psi_{wl}(t)\|_{\s^2_x(\mathbb{R}^n)}\lesssim 1\quad \text{ for all }0\leq k\leq K.
\eeq
Here, the regularity of the interaction determines the value of $K$. In \cite{SW20221}, we extended these results in \cites{LS2020, LS2021} to lower space dimensions and removed the radial symmetry assumption. Specifically, for localized interactions, we have shown that the following asymptotic decomposition holds for all $n\geq 1$:
\eq
\lim\limits_{t\to \infty}\| \psi(t)-e^{-itH_0}\Omega_{+}^*\psi_0-\psi_{wl}(t)\|_{\s^2_x(\mathbb{R}^n)}=0\quad\text{ for all }n\geq 1,
\eeq
where $\psi_{wl}(t)$ satisfies
\eq
(\psi_{wl}(t), |x|\psi_{wl}(t))_{\s^2_x(\mathbb{R}^n)}\lesssim |t|^{1/2+0},\quad t\geq 1.
\eeq
A smooth cut-off function $F_c$ is used to construct the projection on the space of all scattering states in \cite{SW20221}, given by
\eq
P_{sc}(s):=s\text{-}\lim\limits_{t\to \infty} U_{lin}(s,t+s) F_c(\frac{|x-2tP|}{t^\alpha}\leq 1) U_{lin}(s+t,s)\quad \text{ on }\s^2_x(\mathbb{R}^n), \quad n\geq 3, 
\eeq
for all $\alpha \in (0,1-2/n)$, $s\in \mathbb{R}$, where $U_{lin}(a,b)$ denotes the solution operator to the linear Schr\"odinger equation with a time-dependent potential $V(x,t)$. Here $V(x,t)\equiv \mathcal{N}(|\psi(t)|,|x|,t)$ is determined by the solution $\psi(t)$ to the original nonlinear Schr\"odinger equation.

\begin{remark}
 It is NOT true in general that the above projection of the free channel will apply. It can happen that the solution has a big spread  around the classical trajectory, and then the limit will not exist for the values of $\alpha$ above \cites{LS2015,lindblad2023modified,ifrim2022testing}.   
\end{remark}
Once the above construction is done, one can prove local decay estimates. For example, when the potential has the form
\eq
\mathcal{N}=V_0(x)+\left(\sum\limits_{j=1}^{N_1}V_j(x)\sin(\omega_j t)\right)+\left( \sum\limits_{j=N_1+1}^{N}V_j(x)\cos(\omega_j t)\right),
\eeq
where $N, N_1\in \mathbb{N}$ with $N\geq N_1$, and each $V_j(x)$ is localized with $\langle x\rangle^\sigma V_j(x)\in \s^\infty_x(\mathbb{R}^5)$ and $\sigma>6$, we have shown in \cite{Sof-W5} that
\eq
\int_0^\infty dt\|\langle x\rangle^{-\eta}U(t,0)P_{sc}(0)\psi \|_{\s^2_{x}(\mathbb{R}^{5})}^2\lesssim \| \psi\|_{\s^2_x(\mathbb{R}^5)}^2,\text{ for all $\eta>5/2$. }
\eeq

It is worth mentioning that complete results in the radial case, for the NLWE in the scale-invariant energy-critical cases, have been established by C. Kenig, F. Merle, and collaborators in 3 or higher dimensions. See e.g. \cite{duyckaerts2012profiles}.  For further development, see \cite{jendrej2021soliton}.\par

In this paper, we prove the conjecture that $\psi_{wl}(t)$ is localized in $x$ for system \eqref{SE} by improving the decomposition \eqref{eq1} to
\eq
\lim\limits_{t\to \infty}\| \psi(t)-e^{-itH_0}\Omega_{+}^*\psi_0-\psi_{loc}(t)\|_{\s^2_x(\mathbb{R}^n)}=0\label{eq2}
\eeq
for some localized function $\psi_{loc}(t)$ satisfying
\eq
\sup\limits_{t\geq 0}\|\langle x\rangle^\delta \psi_{loc}(t)\|_{\s^2_x(\mathbb{R}^n)}\lesssim_\delta \| \psi_0\|_{\s^2_x(\mathbb{R}^n)}\label{eq3}
\eeq
for some $\delta>0$. One advantage of our method is that we are able to eliminate the component with a frequency of order $\frac{1}{(t+1)^\epsilon}$ for any $\epsilon\in (0,1/2]$. \par
Let $F$ be a smooth characteristic function(a cut-off function) with $F(k)=1$ when $k\geq1$ and $F(k)=0$ when $k\leq 1/2$. Throughout the paper, $F(k>b):=F(\frac{k}{b}),$ $F(k\leq b):=1-F(k>b)$ and $F(b<k\leq c):=F(k>b)-F(k>c)$ for $0<b<c$. \par
Here is {\bf a sketch of the estimate for $F(|P|\leq \frac{1}{(t+1)^\epsilon})\psi(t)$}. Assuming that $V(x,t)$ is well-behaved, we can estimate $F(|P|\leq \frac{1}{(t+1)^\epsilon})\psi(t)$ as follows:
\begin{enumerate}
\item Use incoming/outgoing decomposition:
\eq
F(|P|\leq \frac{1}{(t+1)^\epsilon})\psi(t)=P^+F(|P|\leq \frac{1}{(t+1)^\epsilon})\psi(t)+P^-F(|P|\leq \frac{1}{(t+1)^\epsilon})\psi(t).
\eeq
Here, the definitions of $P^\pm$ can be found in \eqref{Pp} and \eqref{Pm}.
\item For $P^+F(|P|\leq \frac{1}{(t+1)^\epsilon})\psi(t)$, we approximate it using  $P^+F(|P|\leq \frac{1}{(t+1)^\epsilon})\Omega_{t,+}^*\psi(t).$ Similarly, for $P^-F(|P|\leq \frac{1}{(t+1)^\epsilon})\psi(t)$, we use $P^-F(|P|\leq \frac{1}{(t+1)^\epsilon})\Omega_{t,-}^*\psi(t)$, where
\eq
\Omega_{t,\pm}^*\psi(t):=w\text{-}\lim\limits_{s\to \pm\infty} e^{isH_0}\psi(t+s)\quad \text{ on }\s^2_x(\mathbb{R}^n)
\eeq
and $\Omega_{\pm}^*:=\Omega_{0,\pm}^*$, as shown in \eqref{10eq1} and Lemma \ref{Lem2}. Using the intertwining property \eqref{inter} (see Lemma \ref{lem:inter}), we can write:
 \begin{multline}
 P^+F(|P|\leq \\\frac{1}{(t+1)^\epsilon})\psi(t)=P^+F(|P|\leq \frac{1}{(t+1)^\epsilon})\Omega_{t,+}^*\psi(t)+P^+F(|P|\leq \frac{1}{(t+1)^\epsilon})(1-\Omega_{t,+}^*)\psi(t)\\
 =P^+F(|P|\leq \frac{1}{(t+1)^\epsilon})e^{-itH_0}\Omega_+^*\psi(0)+P^+F(|P|\leq \frac{1}{(t+1)^\epsilon})(1-\Omega_{t,+}^*)\psi(t)\\
 =:\psi_1(t)+\psi_2(t).
 \end{multline}
 We then have the following:
 \begin{enumerate}
 \item $\psi_1(t)\to 0$ in $\s^2_x(\mathbb{R}^n)$ as $t\to \infty$ since $ F(|P|\leq \frac{1}{(t+1)^\epsilon})\Omega_+^*\psi(0)\to 0$ in $\s^2_x(\mathbb{R}^n)$ as $t\to \infty$.
 \item For $\psi_2(t)$, we use Duhamel's formula to rewrite it as:
 \eq
 \psi_2(t)=i\int_0^\infty ds P^+F(|P|\leq \frac{1}{(t+1)^\epsilon})e^{isH_0}V(x,t+s)\psi(t+s).
 \eeq
As $t\to \infty$, $\psi_2(t)\to 0$ in $\s^2_x(\mathbb{R}^n)$ for the following two reasons. First, when $n\geq 5$, we have the inequality
\begin{multline}
\int_0^{(t+1)^{2\epsilon}}ds\| e^{isH_0}F(|P|\leq \frac{1}{(t+1)^{\epsilon}}) f(x,s)\|_{\s^2_x(\mathbb{R}^n)}= \\
c\int_0^{(t+1)^{2\epsilon}}ds\| e^{isk^2}F(|k|\leq \frac{1}{(t+1)^{\epsilon}}) \hat{f}(k,s)\|_{\s^2_k(\mathbb{R}^n)}\\
\leq c(t+1)^{2\epsilon}\| F(|k|\leq \frac{1}{(t+1)^{\epsilon}})\|_{\s^2_k(\mathbb{R}^n)}\| \hat{f}(k,t)\|_{\s^\infty_{k,t}(\mathbb{R}^{n+1})}\\
\leq c (t+1)^{2\epsilon}\left( \int d^nk F(|k|\leq \frac{1}{(t+1)^{\epsilon}})^2\right)^{1/2}\| f(x,t)\|_{\s^\infty_t\s^1_{x}(\mathbb{R}^{n+1})}\\
(\text{use }n\geq 5)\leq C \frac{1}{(t+1)^{-2\epsilon+ n/2\times \epsilon}}\| f(x,t)\|_{\s^\infty_t\s^1_{x}(\mathbb{R}^{n+1})}\to 0\label{outline:eq4}
\end{multline}
as $t\to \infty$. Secondly, the inequality
\eq
\int_0^\infty ds \| P^+e^{isH_0}F(|P|\leq 1)\|_{\s^2_{x,\sigma}(\mathbb{R}^n)\cap \s^1_x(\mathbb{R}^n)\to\s^2_x(\mathbb{R}^n) }\lesssim  1,\label{outline:eq1}
\eeq
which holds for all $\sigma>2$, implies that 
\eq
\int_{(t+1)^{2\epsilon}}^\infty ds \| P^+e^{isH_0}F(|P|\leq \frac{1}{(1+t)^\epsilon})\|_{\s^2_{x,\sigma}(\mathbb{R}^n)\cap \s^1_x(\mathbb{R}^n)\to\s^2_x(\mathbb{R}^n) }\to 0
\eeq
as $t\to \infty$. To prove \eqref{outline:eq1}, we establish a low frequency estimate and a high frequency estimate:
\begin{itemize}
\item Low frequency estimate: for all $f(x,t)\in \s^\infty_t\s^1_x(\mathbb{R}^{n+1})$, using the Plancherel theorem and H\"older's inequality, when $n\geq 5$, we have
\begin{multline}
\| e^{isH_0}F(|P|\leq \frac{1}{(s+1)^{1/2-\epsilon}}) f(x,s)\|_{\s^2_x(\mathbb{R}^n)}= \\
c\| e^{isk^2}F(|k|\leq \frac{1}{(s+1)^{1/2-\epsilon}}) \hat{f}(k,s)\|_{\s^2_k(\mathbb{R}^n)}\\
\leq c\| F(|k|\leq \frac{1}{(s+1)^{1/2-\epsilon}})\|_{\s^2_k(\mathbb{R}^n)}\| \hat{f}(k,t)\|_{\s^\infty_{k,t}(\mathbb{R}^{n+1})}\\
\leq c \left( \int d^nk F(|k|\leq \frac{1}{(s+1)^{1/2-\epsilon}})^2\right)^{1/2}\| f(x,t)\|_{\s^\infty_t\s^1_{x}(\mathbb{R}^{n+1})}\\
(\text{use }n\geq 5)\leq C \frac{1}{(s+1)^{n/4- n/2\times \epsilon}}\| f(x,t)\|_{\s^\infty_t\s^1_{x}(\mathbb{R}^{n+1})}\in \s^1_s[0,\infty),\label{outline:eq3}
\end{multline}
for some $c,C>0$. Here, $\hat{f}$ denotes the Fourier transform of $f$ in the $x$ variable. 

\item High frequency estimate: We use a dilation transformation, which keeps the $\s^2_x$ norm invariant, and obtain the following estimate:
\begin{multline}
\| P^+e^{isH_0}F(|P|>\frac{1}{(s+1)^{1/2-\epsilon}})\langle x\rangle^{-\sigma}\|_{\s^2_x(\mathbb{R}^n)\to \s^2_x(\mathbb{R}^n)}=\\
\| P^+e^{i\frac{s}{(1+s)^{1-2\epsilon}}H_0}F(|P|>1)\langle (1+s)^{1/2-\epsilon}x\rangle^{-\sigma}\|_{\s^2_x(\mathbb{R}^n)\to \s^2_x(\mathbb{R}^n)}\\
\leq \|P^+e^{i\frac{s}{(1+s)^{1-2\epsilon}}H_0}F(|P|>1)\chi(|x|\leq 1) \|_{\s^2_x(\mathbb{R}^n)\to\s^2_x(\mathbb{R}^n) }+\\
\| P^+e^{i\frac{s}{(1+s)^{1-2\epsilon}}H_0}F(|P|>1)\langle x\rangle^{-\sigma}\|_{\s^2_x(\mathbb{R}^n)\to \s^2_x(\mathbb{R}^n)}\times\\
\| \chi(|x|\geq 1)\langle x\rangle^\sigma\langle (1+s)^{1/2-\epsilon}x\rangle^{-\sigma}\|_{\s^2_x(\mathbb{R}^n)\to \s^2_x(\mathbb{R}^n)}\in \s^1_s[0,\infty).\label{outline:eq2}
\end{multline}
We obtain \eqref{outline:eq2} by utilizing 
\eq
\|\langle x\rangle^\sigma\langle (1+s)^{1/2-\epsilon}x\rangle^{-\sigma}\|_{\s^2_x(\mathbb{R}^n)\to \s^2_x(\mathbb{R}^n)}\lesssim_n \frac{1}{(1+s)^{(1/2-\epsilon)\sigma}},
\eeq
as well as the estimate \eqref{Sep29.1} from Lemma \ref{out/in1}: for $a\geq 1$ and $\delta>2$,
\eq
\| P^+ e^{ia H_0} \langle x\rangle^{-\delta}\|_{\s^2_x(\mathbb{R}^n)\to\s^2_x(\mathbb{R}^n) }\lesssim_{n,\delta} \frac{1}{a^\delta},
\eeq
\eq
\| P^+ e^{ia H_0} \chi(|x|\leq 1)\|_{\s^2_x(\mathbb{R}^n)\to\s^2_x(\mathbb{R}^n) }\lesssim_{n,N} \frac{1}{a^N}.
\eeq

\end{itemize}
\end{enumerate}
\item Therefore, we have that $ P^+F(|P|\leq \frac{1}{(t+1)^\epsilon})\psi(t)\to 0$ in $\s^2_x(\mathbb{R}^n)$ as $t\to \infty$, and similarly, $P^-F(|P|\leq \frac{1}{(t+1)^\epsilon})\psi(t)\to 0$ in $\s^2_x(\mathbb{R}^n)$ as $t\to \infty$. As a result, we conclude that $F(|P|\leq \frac{1}{(t+1)^\epsilon}) \psi(t)\to 0$ in $\s^2_x(\mathbb{R}^n)$ as $t\to \infty$.
\end{enumerate}
Here is the {\bf detailed outline} of the proof for Theorem \ref{thm} and Theorem \ref{thm2}. The proof of Theorem \ref{thm} is based on the approach initiated in \cite{Sof-W5}. We begin by decomposing the solution into two parts:
\eq
\psi(t)=F(|x|\geq 10)\psi(t)+F(|x|< 10)\psi(t).
\eeq
Note that $F(|x|< 10)\psi(t)$ is localized in $x$. For $F(|x|\geq 10)\psi(t)$, we use an incoming/outgoing decomposition:
\eq
F(|x|\geq 10)\psi(t)=P^+F(|x|\geq 10)\psi(t)+P^-F(|x|\geq 10)\psi(t).
\eeq
To approximate $P^\pm F(|x|\geq 10)\psi(t)$, we use $P^\pm \Omega_{t,\pm}^*\psi(t)$:
\begin{enumerate}
\item A key observation here is that
\eq
\Omega_{t,\pm}^*\psi(t)=w\text{-}\lim\limits_{s\to \pm\infty} e^{isH_0}F(|x|\geq 10)\psi(t+s),
\eeq
as shown in \eqref{10eq1}. Using the intertwining property \eqref{inter}, one has
\eq
P^\pm \Omega_{t,\pm}^*\psi(t)=P^\pm e^{-itH_0}\Omega_{\pm}^*\psi(0),
\eeq
which are close to the free flow. See section \ref{WO} for detailed information about $\Omega_{t,\pm}^*$.
\item We write $P^\pm F(|x|\geq 10)\psi(t)$ as
\eq
P^\pm F(|x|\geq 10)\psi(t)=P^\pm e^{-itH_0}\Omega_{\pm}^*\psi(0)+C_\pm(t)\psi(t)
\eeq
where 
\eq
\begin{split}
C_\pm (t):=&P^\pm F(|x|\geq 10)-P^\pm \Omega_{t,\pm}^*\\
=&\pm i\int_0^\infty ds P^\pm e^{\pm isH_0} F(|x|\geq 10)V(x,t\pm s)U(t\pm s,t)+\\
&(\mp i)\int_0^\infty ds P^\pm e^{\pm is H_0} [H_0, F(|x|\geq 10)] U(t\pm s,t).
\end{split}
\eeq
\item Let $\psi_{free}=\Omega_+^*\psi(0)$ and define $\psi_{loc}(t)$ as $\bar{F}(|x|<10) \psi(t)+(C_+(t)+C_-(t))\psi(t)$. To achieve \eqref{goal}, we need to show that 
\eq
\| \langle x\rangle^\delta(C_+(t)+C_-(t))\psi(t) \|_{\s^2_x(\mathbb{R}^n)}\lesssim_E 1 \text{ for some }\delta>0.
\eeq
This is mainly accomplished by proving that
\eq
\int_0^\infty ds\| \langle x\rangle^\delta P^\pm e^{\pm i sH_0} \|_{\s^2_{x,\sigma}(\mathbb{R}^n)\cap \s^1_x(\mathbb{R}^n)\to\s^2_x(\mathbb{R}^n) }\lesssim_{\sigma,\delta }1\label{4eq1}
\eeq
for some $\delta >0$, where $P:=-i\nabla_x$.
\end{enumerate}
\begin{remark}Here we use the weak limit for $\Omega_{t,\pm}^*$, since the strong limit of $e^{isH_0}\psi(t+s) $ does not exist in general when there is a soliton.
\end{remark}
\begin{remark}Fortunately, the limit $s\text{-}\lim\limits_{s\to \pm\infty} P^\pm e^{isH_0}\psi(t+s) $ exists in $\s^2_x(\mathbb{R}^n)$, therefore there is no confusion in using $\Omega_{t,\pm}^*$.
\end{remark}
\begin{remark}Observe that the weakly localized part only spreads significantly when $|P|$ is close to zero. Here, $P$ refers to the momentum of the weakly localized part. The most challenging part of the argument is showing that the zero-frequency part does not delocalize, that is,
\eq
\| \langle x\rangle^\delta F(|P|\leq \epsilon)\psi_{wl}(t)\|_{\s^2_x(\mathbb{R}^n)}\lesssim 1
\eeq
for some $\delta>0,\epsilon>0$.
\end{remark}
\begin{remark}In $\int_0^\infty ds \| P^\pm e^{\pm isH_0}\|_{\s^2{x,\sigma}(\mathbb{R}^n)\cap \s^1_x(\mathbb{R}^n)\to \s^2_x(\mathbb{R}^n)}$, the part of the solution with momentum (frequency) of order $\frac{1}{\langle s\rangle^{\epsilon}}$ has a total traveling distance $|x|$ of order $\langle s\rangle^{1-\epsilon}\sim s\times |P|$. Using $P^\pm$ and the method of stationary phase, each integration by parts gains $\frac{1}{\langle s\rangle^{1-\epsilon}}$ decay and loses $\langle s\rangle^\epsilon$ due to the cut-off frequency. Therefore, $\epsilon =\frac{1}{2}$ is the borderline. Fortunately, in $5$ or higher space dimensions ($n\geq 5$), we have 
\eq
\| P^\pm e^{\pm isH_0}F(|P|\leq \frac{1}{\langle s\rangle^{1/2-0}})\|_{\s^2_{x,\sigma}(\mathbb{R}^n)\cap \s^1_x(\mathbb{R}^n)\to \s^2_x(\mathbb{R}^n)}\lesssim \frac{1}{\langle s\rangle^{\frac{n}{4}-n/2\times 0 }}\in \s^1_s(\mathbb{R})
\eeq
by using Hardy-Littlewood-Sobolev inequality. So in this case, $n=4$ is the borderline, which is why this method is not workable when $n\leq 4$.
\end{remark}
\begin{remark}
In general, we cannot use propagation estimates with $P^\pm=P^\pm(A)$ (where $A$ is the dilation operator) when $|P|\sim \frac{1}{\langle s\rangle^\epsilon}$, since in this region we have $|A|\lesssim |P||x|\sim \langle s\rangle^{-\epsilon }\cdot \langle s\rangle^\epsilon =O(1)$.
\end{remark}
\begin{remark}Instead of \eqref{H1}, we can use
\eq
\sup\limits_{t\in [0,\infty)}\|\psi(t)\|_{H^1_x(\mathbb{R}^n)}<\infty.\label{H12}
\eeq
we can approximate $P^-F(|x|>10)\psi(t)$ by $P^-F(|x|>10)e^{-itH_0}\psi(0)$, and the error term becomes 
\eq
c\int_0^tdsP^-e^{-i(t-s)H_0}V(x,s)\psi(s)
\eeq
for some $c>0$, using Duhamel's formula. Since $t-s \geq 0$, we can use a similar argument based on estimates of the free flow, and the result follows.
\end{remark}
{The proof of Theorem \ref{thm2} is similar. We use incoming/outgoing decomposition of $\psi(t)$:
\begin{align}
\psi(t)=&P^+\psi(t)+P^-\psi(t)\\
=&P^+e^{-itH_0}\Omega_+^*\psi(0)+P^-e^{-itH_0}\Omega_-^*\psi(0)+P^+(1-\Omega_{t,+}^*)\psi(t)+P^-(1-\Omega_{t,-}^*)\psi(t)\\
=:&P^+e^{-itH_0}\Omega_+^*\psi(0)+P^-e^{-itH_0}\Omega_-^*\psi(0)+C_r(t)\psi(0).
\end{align}
The localized part is defined by
\eq
\psi_{loc}(t)=\tilde{C}_+(t)\psi(0)+\tilde{C}_-(t)\psi(0),
\eeq
where
\eq
\tilde{C}_+(t)\psi(0):=i\int_0^\infty ds P^+e^{isH_0}V_D(x,t+s)\psi_D(t+s)
\eeq
and
\eq
\tilde{C}_-(t)\psi(0):=(-i)\int_0^\infty ds P^-e^{-isH_0}V_D(x,t-s)\psi_D(t-s).
\eeq
We use $\psi_D(t)$ instead of $\psi(t)$ because $\psi_D(t)$ is smoother than $\psi(t)$ when we localize in space (see Lemma \ref{LempsiD}). The difference between $\psi(t)$ and $\psi_D(t)$ is the free flow $e^{-itH_0}\psi(0)$, which is easy to control. In fact, $\langle x\rangle^{-2}\langle P\rangle^{1+l}\psi_D\in\s^2_x(\mathbb{R}^5)$ for any $l\in [0,1)$. In other words, the advantages of using $\psi_D(t)$ to define $\psi_{loc}(t)$ are:
\begin{itemize}
\item $\psi_D(t)$ is smoother than $\psi(t)$.
\item It is easy to control $\psi(t)-\psi_D(t)=e^{-itH_0}\psi(0)$, which satisfies some dispersive estimates (e.g., Strichartz estimate).
\end{itemize}
For $\psi_{loc}(t)$, we obtain \eqref{local} and \eqref{smooth} by using the estimates for $P^\pm e^{\pm isH_0}$ (for $s\geq 0$) acting on localized functions. We are able to prove \eqref{smooth} because $\psi_D(t)$ is smooth.
\begin{remark}Using Duhamel's formula (see \cite{T2006}), we can rewrite $\psi_D(t)$ as follows:
\eq
\psi_D(t)=(-i)\int_0^t ds e^{-i(t-s)H_0}V(x,s)\psi(s).
\eeq
The integration over $s$ makes $\psi_D(t)$ smoother than $\psi(t)$. One can use Duhamel's formula to iterate it again and achieve even greater smoothness. In other words, the power of $2$ in \eqref{smooth} is not optimal, and in most situations, it can be improved to $k>2$ by utilizing the Duhamel iteration strategy.
\end{remark}
}

\section{Preliminaries}\label{sec:Pre}
\subsection{Setup}
Using \eqref{H1}, we can regard the interaction $\mathcal{N}(|\psi|, |x|,t)$ as a general time-dependent perturbation $V(x,t)$, that is,
\eq
V(x,t):=\mathcal{N}(|\psi(t)|, |x|,t).
\eeq
Let us define
\eq
\psi_D(t):=\psi(t)-e^{-itH_0}\psi(0)
\eeq
and
\eq
V_D(x,t):=\mathcal{N}(|\psi_D(t)|,|x|,t).
\eeq
Then, $\psi(t)$ is a solution to the system
\eq
\begin{cases}
    i\partial_t\psi(x,t)=(H_0+V(x,t))\psi(x,t)\\
    \psi(x,0)=\psi_0(x)\in H^1_x(\mathbb{R}^n)
\end{cases},\quad (x,t)\in \mathbb{R}^{n+1}\label{SE0}
\eeq
where $V(x,t)\in \s^\infty_t\s^2_x(\mathbb{R}^{n+1})$ satisfies the following conditions:
\begin{itemize}
\item In Theorem \ref{thm}, \eq
\chi(|x|\geq 1)V(x,t)\psi\in \s^\infty_t\s^2_{x,\sigma}(\mathbb{R}^{n+1})\cap \s^\infty_t\s^1_{x}(\mathbb{R}^{n+1})\label{Vcon}
\eeq
for some $\sigma>2$.
\item In Theorem \ref{thm2},
\eq
|V(x,t)|\lesssim_{E} \frac{1}{\langle x\rangle^{\sigma}}\quad \text{ for some }\sigma >6,
\eeq
\eq
V(x,t)\in L^\infty_tL^2_{x,2}(\R^{n+1}),
\eeq
\begin{multline}
\|V(x,t)\psi(t)-V_D(x,t)\psi_D(t) \|_{\s^2_{x,\sigma}(\mathbb{R}^n)\cap \s^1_{x,2}(\mathbb{R}^n)}\lesssim_{E} \\\|e^{-itH_0}\psi(0)\|_{\s^2_{x,-\sigma'}(\mathbb{R}^n)}+ \|e^{-itH_0}\psi(0)\|_{\s^2_{x,-\sigma'}(\mathbb{R}^n)}^{2/n}\label{Jan26.1}
\end{multline}
for some $\sigma'>0$ and $\sigma>2$,
\eq
\| \langle P\rangle^{3/2}\langle x\rangle^6 V_D(x,t)\psi_D(t)\|_{\s^2_x(\mathbb{R}^n)}\lesssim_{E}1\label{con50}
\eeq
and for all $j=1,2,\cdots,n$,
\eq
\| P_j V(x,t)\psi(t)\|_{\s^\infty_t\s^2_{x,2}(\mathbb{R}^{n+1})}\lesssim_E 1.
\eeq

\end{itemize}

\subsection{Wave operators}\label{WO} The wave operator acting on the initial data $\psi(t)$ for $t\in \mathbb{R}$ is defined as follows:
\eq
Y_{t,\pm}^*\psi(t):=\lim\limits_{s\to \pm \infty} e^{isH_0} U_{lin}(s+t,t)P_c(t)\psi(t)\quad \text{ in }\s^2_x(\mathbb{R}^n).
\eeq
Moreover, we define
\eq
Y_{t,\alpha}^*:=s\text{-}\lim\limits_{s\to \infty} e^{isH_0}F_c(\frac{|x-2sP|}{s^\alpha}\leq 1)U_{lin}(t+s,t)\quad \text{ on }\s^2_x(\mathbb{R}^n),
\eeq
\eq
Y_{t,\alpha}:=s\text{-}\lim\limits_{s\to \infty} U_{lin}(t,t+s)F_c(\frac{|x-2sP|}{s^\alpha}\leq 1)e^{-isH_0}\quad \text{ on }\s^2_x(\mathbb{R}^n)
\eeq
and
\eq
P_{sc}(t):=Y_{t,\alpha}Y_{t,\alpha}^*\quad \text{ on }\s^2_x(\mathbb{R}^n)
\eeq
for $\alpha\in (0, 1-2/n), n\geq3$. Here, $F_c$ denotes a smooth characteristic function, a cut-off function.
\begin{remark}The operator $U_{lin}(t,0)$ is defined as the solution to the Schr\"odinger equation with the potential $V(x,t)$. Hence, the well-definedness of $Y_{t,\alpha}$ is contingent upon determining $V$ through $\psi(t)$. {color{red} We define $Y_{t,\alpha}$ and $Y_{t,\alpha}^*$ solely to give meaning to $P_{sc}$.}
\end{remark}
\begin{lemma}\label{exist:Y}Provided that $V(x,t)\in \s^\infty_t\s^2_x(\mathbb{R}^{n+1})$, for all $\alpha\in (0, 1-2/n)$, $n\geq 3$, $Y_{t,\alpha}$ and $Y_{t,\alpha}^*$ exist on $\s^2_x(\mathbb{R}^n)$ for all $t\in \mathbb{R}$.
\end{lemma}
\begin{proof}The existence of $\Omega_{t,\alpha}$ can be established by employing Cook's method (see, for instance, \cite{S2018}). To obtain the existence of $Y_{t,\alpha}^*$, begin by utilizing Cook's method to express $e^{isH_0}F_c(\frac{|x-2sP|}{s^\alpha}\leq 1)U_{lin}(t+s,t)$ as:
\begin{align}
e^{isH_0}&F_c(\frac{|x-2sP|}{s^\alpha}\leq 1)U_{lin}(t+s,t)=F_c(\frac{|x|}{s^\alpha}\leq 1)e^{isH_0}U_{lin}(t+s,t)\nonumber\\
=&F_c(\frac{|x|}{s^\alpha}\leq 1)e^{isH_0}U_{lin}(t+s,t)\vert_{s=1}+(-i)\int_1^sdu F_c(\frac{|x|}{u^\alpha}\leq 1) e^{iuH_0}V(x,u+t)U_{lin}(t+u,t)+\nonumber\\
&\int_1^sdu \partial_u[F_c(\frac{|x|}{u^\alpha}\leq 1) ]e^{iuH_0}U_{lin}(t+u,t)\nonumber\\
=:&F_c(\frac{|x|}{s^\alpha}\leq 1)e^{isH_0}U_{lin}(t+s,t)\vert_{s=1}+B_{in}(s,t)+B_{\geq0}(s,t).
\end{align}
The interaction term $B_{in}(s,t)$ exists in $\s^2_x(\mathbb{R}^n)$ as $s$ approaches infinity, thanks to the $\s^\infty$ decay estimate of the free flow. As for $B_{\geq0}(s,t)$, the inequality $\partial_u[F_c(\frac{|x|}{u^\alpha}\leq 1) ]\geq 0$ enables us to employ the propagation estimate, which was first introduced in \cite{SW20221}. Using this estimate, we can establish its existence in $\s^2_x(\mathbb{R}^n)$ as $s$ tends to infinity. The key step in proving a propagation estimate is finding a {\emph{Propagation Observable}} (PROB) $B(t)$ bounded(above) s.t.
\eq
\partial_t[(\psi(t), B(t)\psi(t))_{\s^2_x(\mathbb{R}^n)}]=(\psi(t), C^*(t)C(t)\psi(t))_{\s^2_x(\mathbb{R}^n)}+g(t)
\eeq
for some $g(t)\in \s^1_t[1,\infty)$. In our case, $B(t)\equiv F_c(\frac{|x-2tP|}{t^\alpha}\leq 1)$. For further information, please refer to \cite{SW20221}.
\end{proof}
\begin{remark}If we weaken the assumption in Lemma \ref{exist:Y} to $VU_{lin}\phi(x)\in L^\infty_tL^{2n/(n+2)}_x$ for all $\phi(x)\in \s^2_x(\R^n)$, where $n\geq5$, the same argument still holds. See also the Appendix. 
\end{remark}
\begin{lemma}If $V(x,t)\in \s^\infty_t\s^2_x(\mathbb{R}^{n+1})$, then for all $\alpha\in (0, 1-2/n)$, $n\geq 3$, $P_c(t)$ exists on $\s^2_x(\mathbb{R}^n)$ for all $t\in \mathbb{R}$.
\end{lemma}
\begin{proof}The existence of $P_c(t)$ follows from the existence of $Y_{t,\alpha}^*$ and $Y_{t,\alpha}$ on $\s^2_x(\mathbb{R}^n)$.
\end{proof}
In what follows, we shall focus on the nonlinear wave operator $\Omega_{t,\pm}^*$, which acts on the solution $\psi(t)$. $\Omega_{t,\pm}^*\psi(t)$ is defined as 
\eq
\Omega_{t,\pm}^*\psi(t):=s\text{-}\lim\limits_{s\to \pm\infty}e^{isH_0}\psi(t+s),
\eeq
where $U(t,0)$ denotes the solution operator to system \eqref{SE}. {Furthermore, note that $V(x,t)$ satisfies \eqref{Vcon} with $V(x,t) = \mathcal{N}(|\psi(t)|, |x|, t)$, which implies
\eq
\lim\limits_{s\to \infty} P^{+}e^{isH_0} U(s+t,t)\psi(t)\quad \text{ exists in }\s^2_x(\mathbb{R}^n)\label{LemJan1eq}:
\eeq
\begin{lemma}\label{Lem3}If $V(x,t)$ satisfies \eqref{Vcon} and $\psi(t)$ satisfies \eqref{H1}, then \eqref{LemJan1eq} holds for $n\geq 5$.
\end{lemma}
\begin{proof}First, we can observe that:
\eq
\lim\limits_{s\to \infty}  P^{+}e^{isH_0} F(|x|<10)\psi(s+t)=0\quad \text{ in }\s^2_x(\mathbb{R}^n).
\eeq
This follows from \eqref{1eq2} in Lemma \ref{out/in1}, with $l=0$. Therefore, we need only demonstrate the existence of:
\eq
\lim\limits_{s\to \infty} P^{+}e^{isH_0} F(|x|\geq 10)\psi(s+t)\quad \text{ in }\s^2_x(\mathbb{R}^n).
\eeq
Using Cook's method to expand $P^{+}e^{isH_0} F(|x|\geq 10)\psi(s+t)$, we have
\begin{align}
P^{+}e^{isH_0} F(|x|\geq 10)\psi(s+t)=&P^+e^{iH_0}F(|x|\geq 10)\psi(t+1)\nonumber\\
&+(-i)\int_1^s du P^+e^{iuH_0}F(|x|\geq 10) V(x,u+t)\psi(u+t)\nonumber\\
&+i\int_1^s du P^+e^{iuH_0}[H_0,F(|x|\geq 10)] \psi(u+t).
\end{align}
The existence of $\lim\limits_{s\to \infty}P^{+}e^{isH_0} F(|x|\geq 10)\psi(s+t)$ in $\s^2_x(\mathbb{R}^n)$ follows from the fact that, by using \eqref{1eq2} in Lemma \ref{out/in1} and taking $\delta>2$, we obtain
\begin{align}
\| P^+ e^{iuH_0}F(|x|\geq 10)V(x,u+t)\psi(u+t)&\|_{\s^2_x(\mathbb{R}^n)}\leq \| P^+e^{iuH_0}\|_{\s^2_{x,\sigma}(\mathbb{R}^n)\cap\s^1_x(\mathbb{R}^n)\to \s^2_x(\mathbb{R}^n) } \nonumber\\
&\times\|F(|x|\geq 10)V(x,u+t)\psi(t+u)\|_{ \s^2_{x,\sigma}(\mathbb{R}^n)\cap \s^1_x(\mathbb{R}^n)}\nonumber\\
(\text{ use \eqref{Vcon}})\lesssim & \| P^+e^{iuH_0}\|_{\s^2_{x,\sigma}(\mathbb{R}^n)\cap\s^1_x(\mathbb{R}^n)\to \s^2_x(\mathbb{R}^n) }\in \s^1_u[1,\infty)
\end{align}
and
\begin{align}
\| P^+e^{iuH_0}[H_0,F(|x|\geq 10)] \psi(u+t)&\|_{\s^2_x(\mathbb{R}^n)}\leq  \| P^+e^{iuH_0}\|_{\s^2_{x,\sigma}(\mathbb{R}^n)\cap\s^1_x(\mathbb{R}^n)\to \s^2_x(\mathbb{R}^n) } \nonumber\\
&\times\| [H_0, F(|x|\geq 10)]\psi(t+u)\|_{ \s^2_{x,\sigma}(\mathbb{R}^n)\cap\s^1_x(\mathbb{R}^n)}\nonumber\\
\lesssim & \| P^+e^{iuH_0}\|_{\s^2_{x,\sigma}(\mathbb{R}^n)\cap\s^1_x(\mathbb{R}^n)\to \s^2_x(\mathbb{R}^n) }\sup\limits_{t\in \mathbb{R}}\|\psi(t)\|_{H^1_x(\mathbb{R}^n)}\nonumber\\
(\text{ use \eqref{H1}})\lesssim & \| P^+e^{iuH_0}\|_{\s^2_{x,\sigma}(\mathbb{R}^n)\cap\s^1_x(\mathbb{R}^n)\to \s^2_x(\mathbb{R}^n) }\in \s^1_u[1,\infty).
\end{align}
We finish the proof.
\end{proof}
Lemma \ref{Lem3} implies that for $n\geq 5$,
\eq
P^+\Omega_{t,+}^*\psi(t)=P^+\Omega_{t,\alpha}^*\psi(t)=\lim\limits_{s\to \infty} P^{+}e^{isH_0} \psi(s+t)\quad \text{ in }\s^2_x(\mathbb{R}^n).
\eeq
Here, we use $\Omega_{t,\alpha}^*=\Omega_{t,+}^*$ since when $V(x,t)\in \s^\infty_t\s^2_x(\mathbb{R}^{n+1}), n\geq 3$, it was shown in \cite{SW20221} that
\eq
w\text{-}\lim\limits_{s\to \infty} e^{isH_0}(1-F_c(\frac{|x-2sP|}{s^\alpha}\leq 1))U(s+t,t)=0\text{ on }\s^2_x(\mathbb{R}^n).
\eeq}
Similarly, for $n\geq 5$,
\eq
P^-\Omega_{t,-}^*\psi(t)=s\text{-}\lim\limits_{s\to -\infty} P^{-}e^{isH_0} \psi(s+t)\quad \text{ exists in }\s^2_x(\mathbb{R}^n).
\eeq
Let $\Omega_\pm^*:=\Omega_{t,\pm}^*\vert_{t=0}$. The time-dependent wave operators satisfy the following {\bf intertwining property}
\eq
\Omega_{t,\pm}^*U(t,0)=e^{-itH_0}\Omega_\pm^*\label{inter},\text{ on }\s^2_x(\mathbb{R}^n)
\eeq
which follows directly from the definition of $\Omega_{t,\pm}^*$:
\begin{lemma}\label{lem:inter}\eqref{inter} holds if $\Omega_{t,\pm}^*$ exist for all $t\in \mathbb{R}$.
\end{lemma}
\begin{proof}By using the definition of $\Omega_{\pm,t}^*$, we obtain:
\begin{align}
    \Omega_{t,\pm}^*U(t,0)=&s\text{-}\lim\limits_{s\to \pm\infty}e^{isH_0}U(s+t,t)U(t,0)\nonumber\\
    =&e^{-itH_0}s\text{-}\lim\limits_{s\to \pm\infty}e^{i(s+t)H_0}U(s+t,0)\nonumber\\
    =&e^{-itH_0}\Omega_{\pm}^*.
\end{align}
Therefore, the proof is complete.
    
\end{proof}
Furthermore, since
\eq
w\text{-}\lim\limits_{s\to \pm \infty} e^{isH_0}F(|x|<10)U(s+t,t)\psi(t)=0\quad \text{ in }\s^2_x(\mathbb{R}^n),\quad n\geq 1
\eeq
we have
\eq
\Omega_{t,\pm}^*\psi(t)=w\text{-}\lim\limits_{s\to \pm\infty} e^{isH_0}F(|x|\geq 10)\psi(t+s)\quad \text{ in }\s^2_x(\mathbb{R}^n), \quad n\geq 1.\label{10eq1}
\eeq
The following lemma shows that \eqref{10eq1} is indeed true:
\begin{lemma}\label{Lem2}If \eqref{H1}  is satisfied, then \eqref{10eq1} holds.
\end{lemma}
\begin{proof}It follows from
\eq
w\text{-}\lim\limits_{s\to \pm\infty} e^{isH_0}F(|x|<10)\psi(t+s)=0\quad \text{ in }\s^2_x(\mathbb{R}^n), \quad n\geq 1:
\eeq
Choose $f\in \s^2_x(\mathbb{R})$. Using H\"older's inequality and \eqref{H1}, one has
\begin{align}
|(f, e^{isH_0}F(|x|<10)\psi(t+s))_{\s^2_x(\mathbb{R}^n)}|=&|(F(|x|<10)e^{-isH_0}f, \psi(t+s))_{\s^2_x(\mathbb{R}^n)}|\nonumber\\
\leq & \| (F(|x|<10)e^{-isH_0}f\|_{\s^2_x(\mathbb{R}^n)}\|\psi(t+s) \|_{\s^2_x(\mathbb{R}^n)}\nonumber\\
\leq & \| (F(|x|<10)e^{-isH_0}f\|_{\s^2_x(\mathbb{R}^n)}\sup\limits_{t\in \R}\|\psi(t)\|_{H^1_x(\mathbb{R}^n)}\nonumber\\
\to & 0,
\end{align}
as $s\to \pm \infty$.
\end{proof}

\subsection{Incoming/outgoing decomposition}\label{subsec2.3}
The incoming/outgoing wave decompositions are similar to the ones initiated by Mourre \cite{M1979}. Let $A$ be the dilation operator
\eq
A:=\frac{1}{2}(x\cdot P+P\cdot x).
\eeq
\begin{definition}[Incoming/outgoing waves] The projection on outgoing waves is defined by (see \cite{S2011})
\eq
P^+:= (\tanh(\frac{A-M}{R})+1)/2\label{Pp}
\eeq
for some sufficiently large $R,M>0$ such that Lemma \ref{PPpm} holds. The projection on incoming waves is defined by
\eq
P^-:=1-P^+.\label{Pm}
\eeq
\end{definition}
Recall that $F$ is a smooth characteristic function(a cut-off function) with $F(k)=1$ when $k\geq1$ and $F(k)=0$ when $k\leq 1/2$. Throughout the paper, $F(k>b):=F(\frac{k}{b}),$ $F(k\leq b):=1-F(k>b)$ and $F(b<k\leq c):=F(k>b)-F(k>c)$ for $0<b<c$.
\begin{lemma}\label{PPpm}Given $\sigma>\max\{1,n/2\}$, for all $N\in [1,\sigma]$, $R>2N/\pi$,  $f\in \s^2_x(\mathbb{R}^n)$, and $n\geq 1$,
\eq
\| \langle x\rangle^N P^\pm \langle x\rangle^{-N}  f\|_{\s^2_x(\mathbb{R}^n)}\lesssim_{\pi/2-N/R,N/R}\|  f\|_{\s^2_x(\mathbb{R}^n)}.\label{PPpmeq}
\eeq
\end{lemma}
\begin{proof}First, we have
\begin{align}
\| \langle x\rangle^N P^\pm \langle x\rangle^{-N}  f\|_{\s^2_x(\mathbb{R}^n)}\leq &\|\chi(|x|\leq 1)  \langle x\rangle^N P^\pm \langle x\rangle^{-N}  f\|_{\s^2_x(\mathbb{R}^n)}+\|\chi(|x|> 1)  \langle x\rangle^N P^\pm \langle x\rangle^{-N}  f\|_{\s^2_x(\mathbb{R}^n)}\nonumber\\
\lesssim & \| P^\pm \langle x\rangle^{-N}  f\|_{\s^2_x(\mathbb{R}^n)}+\|\chi(|x|> 1)  \langle x\rangle^N P^\pm \langle x\rangle^{-N}  f\|_{\s^2_x(\mathbb{R}^n)}\nonumber\\
\lesssim & \| f\|_{\s^2_x(\mathbb{R}^n)}+\|\chi(|x|> 1)  \langle x\rangle^N P^\pm \langle x\rangle^{-N}  f\|_{\s^2_x(\mathbb{R}^n)}.\label{Jan300}
\end{align}
{\bf Estimate for $\|\chi(|x|> 1)  \langle x\rangle^N P^\pm \langle x\rangle^{-N}  f\|_{\s^2_x(\mathbb{R}^n)}$:} We compute $[|x|^N, P^+] \langle x\rangle^{-N}f$ as follows:
\begin{align}
[|x|^N, P^+] \langle x\rangle^{-N}  f=&c\int dw \hat{g}(w)( |x|^N e^{i\omega(\frac{A-M}{R})} -e^{i\omega(\frac{A-M}{R})} |x|^N)\langle x\rangle^{-N}f(x)\nonumber\\
=&c\int dw \hat{g}(w) e^{i\frac{w(A-M)}{R}} ( \frac{1}{e^{Nw/R}}-1)|x|^N\langle x\rangle^{-N} f(x)
\end{align}
for some $c>0$. Here $\hat{g}(w)$ denotes the Fourier transform of $g(k)$ in $k$ variable($k\in \mathbb{R}$).
Since
\eq
\hat{g}(w)=\frac{c_1}{\sinh(\pi w/2)}+c_2\delta(w)\quad \text{ for some constants }c_1,c_2,
\eeq
we have
\begin{align}
\|[|x|, P^+] \langle x\rangle^{-1}  f\|_{\s^2_x(\mathbb{R}^n) }\lesssim &2cc_1\int_{|w|\geq 1} dw | \frac{1}{\sinh(\pi w/2)}| e^{N|w|/R}\| f\|_{\s^2_x(\mathbb{R}^n)}+\\
&cc_1\int_{|w|< 1} dw |\hat{g}(w) || ( \frac{1}{e^{Nw/R}}-1)|\| f\|_{\s^2_x(\mathbb{R}^n)}\nonumber\\
\lesssim_{N/R}& \int_{|w|\geq 1} dw | \frac{1}{\sinh(\pi w/2)}| e^{N|w|/R}\| f\|_{\s^2_x(\mathbb{R}^n)}+\int_{|w|< 1} dw |w\hat{g}(w) |\| f\|_{\s^2_x(\mathbb{R}^n)}\nonumber\\
\lesssim_{\pi/2-N/R,N/R}& \| f\|_{\s^2_x(\mathbb{R}^n)},
\end{align}
when $R\pi/2>1$. Thus,
\eq
\|[|x|^N, P^-] \langle x\rangle^{-N}  f\|_{\s^2_x(\mathbb{R}^n) }=\|[|x|^N, P^+] \langle x\rangle^{-N}  f\|_{\s^2_x(\mathbb{R}^n) }\lesssim_{\pi/2-N/R,R} \| f\|_{\s^2_x(\mathbb{R}^n)},
\eeq
which implies
\eq
\begin{split}
\| \chi(|x|>1)\langle x\rangle^N P^+ \langle x\rangle^{-N}f\|_{\s^2_x(\mathbb{R}^n)}\leq & \| \chi(|x|>1)\frac{\langle x\rangle^N}{|x|^N}[|x|^N, P^+] \langle x\rangle^{-N}f\|_{\s^2_x(\mathbb{R}^n)}+\\
&\| \chi(|x|>1)\frac{\langle x\rangle^N}{|x|^N} P^+\frac{|x|^N}{\langle x\rangle^N}f \|_{\s^2_x(\mathbb{R}^n)}\\
\lesssim_{\pi/2-N/R,R}& \| f\|_{\s^2_x(\mathbb{R}^n)}.
\end{split}\label{Jan301}
\eeq
Similarly, we have
\eq
\| \chi(|x|>1)\langle x\rangle^N P^- \langle x\rangle^{-N}f\|_{\s^2_x(\mathbb{R}^n)}\lesssim_{\pi/2-N/R,R} \| f\|_{\s^2_x(\mathbb{R}^n)}.\label{Jan302}
\eeq
By using \eqref{Jan300}, \eqref{Jan301} and \eqref{Jan302}, we obtain \eqref{PPpmeq} and complete the proof.

\end{proof}
Incoming and outgoing waves exhibit the following properties, which have been shown in \cite{HSS1999} for energies that are away from both $0$ and $\infty$.
\begin{lemma}\label{out/in1}When $R> 2/\pi$, the incoming and outgoing waves satisfy:
\begin{enumerate}
\item {\bf High Energy Estimate}

For all $\sigma> 1, t\geq 1, c>0$, when the space dimension $n\geq 1$, $l\in [0,\sigma)$,
\eq
\|P^{\pm}  F(|P|>c)e^{\pm i tH_0}|P|^l\langle x\rangle^{-\sigma}\|_{\s^2_x(\mathbb{R}^n)\to \s^2_x(\mathbb{R}^n)}\lesssim_{c,n,l,\sigma}  \frac{1}{\langle t\rangle^{\sigma}} .\label{Sep29.1}
\eeq
\item {\bf Near Threshold Estimate}

For $\sigma>1$ when $t\geq1,$ ($\epsilon\in (0,1/2)$), $l\in [0,\sigma)$,
\eq
\| P^{\pm}  F(|P|>\frac{1}{\langle t\rangle^{1/2-\epsilon}})e^{\pm i tH_0} |P|^l\langle x\rangle^{-\sigma}\|_{\s^2_x(\mathbb{R}^n)\to \s^2_x(\mathbb{R}^n)}\lesssim_{n,\epsilon,l}  \frac{1}{t^{(1/2+\epsilon)\sigma+(1/2-\epsilon)l}}  .\label{Sep20.1}
\eeq
In particular, when $n\geq 5, \sigma>2,$ and $ l\in [0,\sigma)$, one has that for some $\epsilon>0$,
\begin{multline}
\| P^{\pm}  e^{\pm i tH_0} |P|^l \|_{\s^2_{x,\sigma}(\mathbb{R}^n)\cap \s^1_x(\mathbb{R}^n)\to \s^2_x(\mathbb{R}^n)}\lesssim_{n,\epsilon,l,\sigma} \\
\max\{\frac{\chi(t\geq 1)}{\langle t\rangle^{(1/2-\epsilon)(n/2+l)}}, \frac{\chi(t\geq 1)}{\langle t\rangle^{\sigma(1/2+\epsilon)+l(1/2-\epsilon)}}\}\in \s^1_t(\mathbb{R})\label{1eq2}.
\end{multline}
This is because the $\s^2 $ volume of $|P|^lF(|P|\leq \frac{1}{\langle t\rangle^{1/2-\epsilon}})f$ in frequency space is controlled by $$\|f\|_{\s^1_x}/\langle t\rangle^{(1/2-\epsilon)(l+n/2)}$$up to some constant.
\item{\bf Time Smoothing Estimate}
For $\sigma>2$, $c>0$, $n\geq 5, l\in [0,1)$, and $a=0,1,2$, we have
\eq
\int_0^\infty  dt \|P^\pm e^{\pm itH_0}|P|^{l} \|_{\s^2_{x,\sigma}(\mathbb{R}^n)\cap \s^1_x(\mathbb{R}^n)\to \s^2_x(\mathbb{R}^n)}\lesssim_{n,l,\sigma} 1\label{1eq3},
\eeq
\eq
\int_0^\infty  dt \|P^\pm e^{\pm itH_0}|P|^{1/2} \|_{\s^2_{x,\sigma}(\mathbb{R}^n)\to \s^2_x(\mathbb{R}^n)}\lesssim_{n,\sigma} 1\label{1eq6}.
\eeq
When $\sigma>3$, $c>0$, $n\geq 5$, $l\in [0,1)$, and $a=0,1,2$, we have
\eq
\int_0^1 dt t^a\|P^{\pm}  F(|P|>c)e^{\pm i tH_0}|P|^{a+l}\|_{\s^2_{x,\sigma}(\mathbb{R}^n)\to \s^2_x(\mathbb{R}^n)}\lesssim_{c,n,l,a} 1 .\label{Oct.1}
\eeq

\item {\bf Weight Absorption Estimate }
\begin{enumerate}
\item For $\sigma>2$, $n\geq5$, $\delta\in(0,   \min\{ \sigma/20-1/10,1/40\}) $, one has
\eq
\int_0^\infty dt\| \langle x\rangle^\delta P^\pm e^{\pm itH_0}\|_{\s^2_{x,\sigma}(\mathbb{R}^n)\cap \s^1_x(\mathbb{R}^n)\to\s^2_x(\mathbb{R}^n) }\lesssim 1.\label{1eq1}
\eeq
Here, it should be pointed out that this estimate is different from the standard propagation estimates done via Mourre's method. In this estimate the thresholds $0,\infty$ are included.
\item For $\sigma>n/2$, $n\geq 45$ and $R>2\sigma/\pi$, there exists $\delta>1$ such that
\eq
\int_0^\infty dt\| \langle x\rangle^\delta P^\pm e^{\pm itH_0}\|_{\s^2_{x,\sigma}(\mathbb{R}^n)\cap \s^1_x(\mathbb{R}^n)\cap H^{3/2}_x(\mathbb{R}^n)\to\s^2_x(\mathbb{R}^n) }\lesssim 1.\label{100eq1}
\eeq
\end{enumerate}
\item {\bf Global Time Smoothing Estimate }For $a=0,1,2,$ $\sigma>4$, $n\geq 5$,
\eq
 \int_0^\infty ds s^a\| P^+e^{isH_0}(-\Delta)^aF(|P|\leq 1)\|_{\s^2_{x,\sigma}(\mathbb{R}^n)\cap \s^1_x(\mathbb{R}^n)\to\s^2_x(\mathbb{R}^n) }\lesssim_{\sigma,n} 1.\label{Jan23.1}
 \eeq
\end{enumerate}

\end{lemma}

\begin{proof}
It suffices to check the case when $t>1$. Let
\eq
(\textbf{Littlewood–Paley Projections})\quad F_{2^j}(|P|):= F(2^j<|P|\leq  2^{j+1}),\quad j=1,2,\cdots.
\eeq
Here $F_{2^j}(|P|)$ satisfies
\eq
\sum\limits_{j=0}^\infty F_{2^j}(|P|)=F(|P|>1).
\eeq
It is known that
\eq
\|P^{\pm}  F_{2^0}(|P|)e^{\pm i tH_0}\langle x\rangle^{-\sigma}\|_{\s^2_x(\mathbb{R}^n)\to \s^2_x(\mathbb{R}^n)}\lesssim_{n}  \frac{1}{\langle t\rangle^{\sigma}},
\eeq
see \cite{HSS1999}. \\
\textbf{Proof of \eqref{Sep29.1}:} When $c=1$, break $P^{\pm}  F(|P|>1)e^{\pm i tH_0}|P|^l\langle x\rangle^{-\sigma}$ into several pieces
\eq
F(|P|>1)e^{\pm i tH_0}|P|^l\langle x\rangle^{-\sigma}=\sum\limits_{j=0}^\infty P^{\pm}  F_{2^j}(|P|)|P|^le^{\pm i tH_0}\langle x\rangle^{-\sigma}.
\eeq
Let $\tilde{F}(k)\in C_0(\mathbb{R})$ be a smooth cut-off function such that $\tilde{F}(k)F_{2^0}(k)=F_{2^0}(k)$. Using dilation transformation, one has
\begin{multline}
\| P^{\pm}  F_{2^j}(|P|)e^{\pm i tH_0}|P|^l\langle x\rangle^{-\sigma}\|_{\s^2_x(\mathbb{R}^n)\to \s^2_x(\mathbb{R}^n)}=\| P^{\pm}  F_{1}(|P|)e^{\pm i 2^{2j}tH_0}|2^jP|^l\langle x/2^j\rangle^{-\sigma}\|_{\s^2_x(\mathbb{R}^n)\to \s^2_x(\mathbb{R}^n)}\\
\leq 2^{jl}\| P^{\pm}  F_{1}(|P|)e^{\pm i 2^{2j}tH_0}\langle x\rangle^{-\sigma}\|_{\s^2_x(\mathbb{R}^n)\to\s^2_x(\mathbb{R}^n) } \times\| \langle x\rangle^\sigma |P|^l\tilde{F}(|P|)\langle x/2^j\rangle^{-\sigma}\|_{\s^2_x(\mathbb{R}^n)\to \s^2_x(\mathbb{R}^n)}\\
(\text{use }t\geq 1)\lesssim_n \frac{1}{\langle 2^{2j}t\rangle^\sigma} \times 2^{j(\sigma+l)}\lesssim_n \frac{1}{2^{j(\sigma-l)}} \times \frac{1}{\langle t\rangle^{\sigma}},
\end{multline}
which implies
\eq
\| \sum\limits_{j=0}^\infty P^{\pm}  F_{2^j}(|P|)e^{\pm i tH_0}|P|^l\langle x\rangle^{-\sigma}\|_{\s^2_x(\mathbb{R}^n)\to \s^2_x(\mathbb{R}^n)}\lesssim_{n,l,\sigma} \frac{1}{\langle t\rangle^\sigma}.
\eeq
We finish the proof for \eqref{Sep29.1} when $c=1$. When $c\neq 1$, \eqref{Sep29.1} will follow in a similar way by using dilation transformation.\\
\textbf{Proof of \eqref{Sep20.1}: }
Take $\tilde{t}:=\left(t/\langle t\rangle^{1-2\epsilon}\right)^{1/(2\epsilon)}$. Then $\tilde{t}\sim t$ when $t\geq 1$. Using dilation transformation, one has
\begin{multline}
\| P^{\pm}  F(|P|>\frac{1}{\langle t\rangle^{1/2-\epsilon}})e^{\pm i tH_0}|P|^l \langle x\rangle^{-\sigma}\|_{\s^2_x(\mathbb{R}^n)\to \s^2_x(\mathbb{R}^n)}=\\
\frac{1}{\langle t\rangle^{(1/2-\epsilon)l}} \| P^{\pm}  F(|P|>1)e^{\pm i \tilde{t}^{2\epsilon}H_0} |P|^l\langle \langle t\rangle^{1/2-\epsilon} x\rangle^{-\sigma}\|_{\s^2_x(\mathbb{R}^n)\to \s^2_x(\mathbb{R}^n)} \lesssim_{n,\epsilon}  \\
\frac{1}{\langle t\rangle^{(1/2-\epsilon)l}} \| P^{\pm}  F(|P|>1)e^{\pm i \tilde{t}^{2\epsilon}H_0} |P|^l\langle \langle t\rangle^{1/2-\epsilon} x\rangle^{-\sigma}\chi(|x|\leq 1)\|_{\s^2_x(\mathbb{R}^n)\to \s^2_x(\mathbb{R}^n)}+\\
\frac{1}{\langle t\rangle^{(1/2-\epsilon)l}} \| P^{\pm}  F(|P|>1)e^{\pm i \tilde{t}^{2\epsilon}H_0} |P|^l\langle \langle t\rangle^{1/2-\epsilon} x\rangle^{-\sigma}\chi(|x|>1)\|_{\s^2_x(\mathbb{R}^n)\to \s^2_x(\mathbb{R}^n)} \\\lesssim_{n,\epsilon}\frac{1}{\langle t\rangle^{(1/2+\epsilon)\sigma+l(1/2-\epsilon)}}.
\end{multline}
Here we use the fact that when $|x|\geq 1$ on the right-hand side, using \eqref{Sep29.1}, $\langle x\rangle^\sigma\langle \langle t\rangle^{1/2-\epsilon}x\rangle^{-\sigma}$ gives $ \frac{1}{\langle t\rangle^{(1/2-\epsilon)\sigma}}$ decay and
\eq
\| P^\pm e^{\pm it^{2\epsilon}H_0}\langle x\rangle^{-\sigma}\|_{\s^2_x(\mathbb{R}^n)\to\s^2_x(\mathbb{R}^n) }
\eeq
gives $1/\langle t\rangle^{2\epsilon \sigma}$ decay. When $|x|<1$, the right-hand side is well-localized and one gets the same decay at least by applying \eqref{Sep29.1}. So we get \eqref{Sep20.1}.\\
\textbf{Proof of \eqref{1eq3}:} Based on\eqref{1eq2}, it suffices to check
\eq
B_1:=\int_0^1  dt \|P^\pm e^{\pm itH_0}|P|^{l}F(|P|\geq 1)\langle x\rangle^{-\sigma} \|_{\s^2_x(\mathbb{R}^n)\to \s^2_x(\mathbb{R}^n)}
\eeq
and
\eq
B_2:=\int_0^1  dt \|P^\pm e^{\pm itH_0}|P|^{l}F(|P|< 1) \|_{\s^1_x(\mathbb{R}^n)\to \s^2_x(\mathbb{R}^n)}.
\eeq
For $B_2$, one has
\eq
B_2\leq  \| |P|^l F(|P|<1)\|_{\s^1_x(\mathbb{R}^n)\to \s^2_x(\mathbb{R}^n)}\lesssim 1.
\eeq
For $B_1$, one has
\begin{align}
B_1\leq & \sum\limits_{j=0}^\infty \int_0^1 dt \| P^\pm e^{\pm i tH_0} |P|^l F_j(|P|)\langle x\rangle^{-\sigma}\|_{\s^2_x(\mathbb{R}^n)\to \s^2_x(\mathbb{R}^n)}\nonumber\\
\leq &\sum\limits_{j=0}^\infty \int_0^1 dt \| P^\pm e^{\pm i tH_0} |P|^l F_j(|P|)\langle x\rangle^{-\sigma}\|_{\s^2_x(\mathbb{R}^n)\to \s^2_x(\mathbb{R}^n)}\chi(t\leq \frac{1}{2^j})\nonumber\\
&+\sum\limits_{j=0}^\infty \int_0^1 dt \| P^\pm e^{\pm i tH_0} |P|^l F_j(|P|)\langle x\rangle^{-\sigma}\|_{\s^2_x(\mathbb{R}^n)\to \s^2_x(\mathbb{R}^n)}\chi(t> \frac{1}{2^j})\nonumber\\
=:&B_{11}+B_{12}.
\end{align}
For $B_{11}$, one has
\begin{align}
B_{11}\lesssim  & \sum\limits_{j=0}^\infty \int_0^1 dt \chi(t\leq \frac{1}{2^j})\times 2^{jl}\nonumber\\
\lesssim&\sum\limits_{j=1}^\infty \frac{1}{2^{j(1-l)}}\nonumber\\
\lesssim_l & 1.\label{B_11}
\end{align}
For $B_{12}$, using dilation transformation and {\bf High Energy Estimate \eqref{Sep29.1}}, one has
\begin{align}
B_{12}= &\sum\limits_{j=0}^\infty \int_0^1 dt \| P^\pm e^{\pm i t2^{2j}H_0} |2^jP|^l F_0(|P|)\langle x/2^j\rangle^{-\sigma}\|_{\s^2_x(\mathbb{R}^n)\to \s^2_x(\mathbb{R}^n)}\chi(t> \frac{1}{2^j})\nonumber\\
\lesssim &\sum\limits_{j=0}^\infty  \int_0^1 dt 2^{jl}\times\frac{1}{\langle 2^{2j}t\rangle^\sigma}\|\langle x\rangle^\sigma \langle x/2^j\rangle^{-\sigma} \|_{\s^2_x(\mathbb{R}^n)\to\s^2_x(\mathbb{R}^n) }\chi(t> \frac{1}{2^j})\nonumber\\
\lesssim &\sum\limits_{j=0}^\infty \int_0^1 \frac{1}{2^{j(\sigma-l)}}\times \frac{1}{t^\sigma}\chi(t> \frac{1}{2^j})\nonumber\\
\lesssim& \sum\limits_{j=0}^\infty \frac{2^{j(\sigma-1)}}{2^{j(\sigma-l)}}\nonumber\\
\lesssim_l & 1.\label{B_12}
\end{align}
Based on \eqref{B_11} and \eqref{B_12}, we have
\eq
B_1\lesssim_l 1.\label{eB1}
\eeq
We finish the proof for \eqref{1eq3}. \\
\textbf{Proof of \eqref{1eq6}:} According to \eqref{eB1}, one has
\eq
\int_0^1dt \|P^\pm e^{\pm i tH_0}|P|^{1/2}\|_{\s^2_{x,\sigma}(\mathbb{R}^n)\to \s^2_x(\mathbb{R}^n)}\lesssim_{n,\sigma} 1.\label{3eq1}
\eeq
For $t\geq 1$, using \eqref{Sep20.1}, we can derive that for $n\geq 5$,
\eq
\| P^{\pm}  e^{\pm i tH_0} |P|^{1/2} \|_{\s^2_{x,\sigma}(\mathbb{R}^n)\to \s^2_x(\mathbb{R}^n)}\lesssim_{n,\epsilon} \max\{\frac{\chi( t\geq 1)}{\langle t\rangle^{5/4-5/2\epsilon}}, \frac{\chi(t\geq 1)}{\langle t\rangle^{\sigma(1/2+\epsilon)+1/2(1/2-\epsilon)}}\}\in \s^1_t(\mathbb{R})\label{3eq2}
\eeq
since for $f\in \s^2_{x,\sigma}(\mathbb{R}^n)$, using the Plancherel theorem and H\"older's inequality, we obtain
\begin{align}
\| |P|^{1/2}F(|P|\leq \frac{1}{\langle t\rangle^{1/2-\epsilon}})f\|_{\s^2_x(\mathbb{R}^n)}=&c_n\| |q|^{1/2}F(|q|\leq \frac{1}{\langle t\rangle^{1/2-\epsilon}})\hat{f}(q)\|_{\s^2_q(\mathbb{R}^n)}\nonumber\\
\lesssim& \frac{1}{\langle t\rangle^{1/4-\epsilon/2}}\|F(|q|\leq \frac{1}{\langle t\rangle^{1/2-\epsilon}}) \|_{\s^{2/n}_q(\mathbb{R}^n)}\| \hat{f}(q)\|_{\s^{(n-4)/2n}_q(\mathbb{R}^n)}\nonumber\\
\lesssim &\frac{1}{\langle t\rangle^{5/4-5/2\epsilon}}\| f(x)\|_{\s^{(n+4)/2n}_x(\mathbb{R}^n)}\nonumber\\
(\text{use }\sigma >2)\lesssim & \frac{1}{\langle t\rangle^{5/4-5/2\epsilon}}\|f(x)\|_{\s^2_{x,\sigma}(\mathbb{R}^n)}.
\end{align}
Here $\hat{f}$ denotes the Fourier transform of $f$ in $x$ variable. Using \eqref{3eq1} and \eqref{3eq2}, one gets \eqref{1eq6} and finish the proof.\\
\textbf{Proof of \eqref{Oct.1}:} When $c=1$, estimate the left-hand side(LHS) of \eqref{Oct.1}
\begin{multline}
\text{ (LHS) of \eqref{Oct.1}}\leq \sum\limits_{j=0}^\infty\int_0^1t^{a}dt\| P^{\pm}  F_{2^j}(|P|)e^{\pm i tH_0}|P|^{a+l}\langle x\rangle^{-\sigma}\|_{\s^2_x(\mathbb{R}^n)\to \s^2_x(\mathbb{R}^n)} \\
\leq\sum\limits_{j=0}^\infty\int_0^{1/2^{j}}t^adt\| P^{\pm}  F_{2^j}(|P|)e^{\pm i tH_0}|P|^{a+l}\langle x\rangle^{-\sigma}\|_{\s^2_x(\mathbb{R}^n)\to \s^2_x(\mathbb{R}^n)}+\\
\sum\limits_{j=0}^\infty\int_{1/2^{j}}^1t^adt\| P^{\pm}  F_{2^j}(|P|)e^{\pm i tH_0}|P|^{a+l}\langle x\rangle^{-\sigma}\|_{\s^2_x(\mathbb{R}^n)\to \s^2_x(\mathbb{R}^n)}=:\sum\limits_{j=0}^\infty A_{j,1}+\sum\limits_{j=0}^\infty A_{j,2}.
\end{multline}
For $A_{j,1}$, one has
\eq
\sum\limits_{j=0}^\infty A_{j,1}\lesssim\sum\limits_{j=0}^\infty  \int_0^{1/2^{j}}dt t^a \times 2^{j(a+l)}\lesssim \sum\limits_{j=0}^\infty 2^{-j(1-l)}\lesssim_l 1.
\eeq
For $A_{j,2}$, we can use dilation to replace $|P|$ with $2^j|P|$. By applying the {\bf High Energy Estimate \eqref{Sep29.1}}, we have:
\begin{align}
A_{j,2}=&2^{(a+l)j}\int_{1/2^{j}}^1 t^adt \| P^{\pm}  F_{1}(|P|)e^{\pm i 2^{2j}tH_0}|P|^{a+l}\langle x/2^j\rangle^{-\sigma}\|_{\s^2_x(\mathbb{R}^n)\to \s^2_x(\mathbb{R}^n)}\\
\lesssim_{n,a,l,\sigma}& 2^{(a+l)j}\int_{1/2^{j}}^1 t^adt\times \frac{1}{\langle t2^{2j}\rangle^\sigma}\times 2^{j\sigma}\\
\lesssim_{n,a,l,\sigma}& \frac{1}{2^{j(1-l)}}+\frac{1}{2^{j(\sigma-a-1)}}.
\end{align}
This expression is summable over $j$ when $\sigma >3$, $a=0,1,2,$ and $l\in [0,1)$. Based on the estimates on $A_{j,1}$ and $A_{j,2}$, one gets \eqref{Oct.1} when $c=1$. For the case with $c>0$, the argument is the same (One can use dilation to reduce the problem to the case $c=1$). \\
\textbf{Proof of \eqref{1eq1}:} When $t\geq 0$, let
\eq
B_3:=\int_0^\infty \|F(|x|\leq 1)\langle x\rangle^\delta P^\pm e^{\pm itH_0}\|_{\s^2_{x,\sigma}(\mathbb{R}^n)\cap \s^1_x(\mathbb{R}^n)\to\s^2_x(\mathbb{R}^n) }dt,
\eeq
\eq
B_4:=\sum\limits_{j=0}^\infty\int_0^\infty \|F_{2^j}(|x|)\langle x\rangle^\delta P^\pm e^{\pm itH_0}\|_{\s^2_{x,\sigma}(\mathbb{R}^n)\cap \s^1_x(\mathbb{R}^n)\to\s^2_x(\mathbb{R}^n) } \chi(t\geq 2^{j/10})dt,
\eeq
and
\eq
B_5:=\sum\limits_{j=0}^\infty\int_0^\infty \|F_{2^j}(|x|)\langle x\rangle^\delta P^\pm e^{\pm itH_0}\langle x\rangle^{-\sigma}\|_{\s^2_x(\mathbb{R}^n)\to\s^2_x(\mathbb{R}^n) } \chi(t<2^{j/10})dt.
\eeq
Then
\eq
\text{(LHS) of \eqref{1eq1}}\leq B_3+B_4+B_5.
\eeq
For $B_3$ and $B_4$, according to \eqref{1eq3} and \eqref{1eq2}, when $n\geq5$, we have
\eq
B_3\lesssim \int_0^\infty dt\| P^\pm e^{\pm itH_0}\|_{\s^2_{x,\sigma}(\mathbb{R}^n)\cap \s^1_x(\mathbb{R}^n)\to\s^2_x(\mathbb{R}^n) }\lesssim1\label{B3}
\eeq
and
\begin{align}
B_4\lesssim &\sum\limits_{j=0}^\infty\int_0^\infty dt \chi(t\geq 2^{j/10})2^{j\delta}\| P^\pm e^{\pm itH_0}\|_{\s^2_{x,\sigma}(\mathbb{R}^n)\cap \s^1_x(\mathbb{R}^n)\to\s^2_x(\mathbb{R}^n) }\nonumber\\
\lesssim_\epsilon&1+ \int_1^\infty t^{10\delta}\times\max\{ \frac{1}{\langle t\rangle^{n/4-n\times\epsilon/2}},\frac{1}{\langle t\rangle^{\sigma(1/2+\epsilon)}}\}\nonumber\\
\lesssim_{\sigma,\delta,\epsilon} & 1\label{B4}
\end{align}
by choosing $\epsilon >0$ small enough such that
\eq
\sigma(1/2+\epsilon)-1>10\delta\label{B41}
\eeq
and
\eq
n/4-n\epsilon/2-1>10\delta.\label{B42}
\eeq
Here, when $n$ is greater than or equal to 5, all we need to do is verify that
\eq
5/4-5\epsilon/2-1>10\delta,
\eeq
and we will also utilize the fact that $\delta\in (0,\min\{\sigma/20-1/10,1/40\})$.\\
For $B_5$,
\begin{align}
B_5\leq& \sum\limits_{j=0}^\infty \int_0^\infty dt \| \langle x\rangle^{\delta} F_{2^j}(|x|)P^\pm e^{i\pm tH_0}F(|P|\leq 2^{j/2}) \langle x\rangle^{-\sigma}\|_{\s^2_x(\mathbb{R}^n)\to\s^2_x(\mathbb{R}^n) } \chi(t<2^{j/10})\nonumber\\
&+\sum\limits_{j=0}^\infty \int_0^\infty dt \| \langle x\rangle^\delta F_{2^j}(|x|)P^\pm e^{i\pm tH_0}F(|P|> 2^{j/2}) \langle x\rangle^{-\sigma}\|_{\s^2_x(\mathbb{R}^n)\to\s^2_x(\mathbb{R}^n) } \chi(t<2^{j/10})\nonumber\\
=:&B_{51}+B_{52}.
\end{align}
For $B_{51}$, when $\delta<1/40$, we can use Lemma \ref{PPpm} to obtain
\begin{align}
B_{51}\lesssim &\sum\limits_{j=0}^\infty \int_0^{2^{j/10}} dt 2^{(\delta-1) j}\| \langle x \rangle P^\pm \langle x\rangle^{-1}\|_{\s^2_x(\mathbb{R}^n)\to\s^2_x(\mathbb{R}^n) }\| \langle x\rangle e^{\pm itH_0}F(|P|\leq 2^{j/2}) \langle x\rangle^{-\sigma}\|_{\s^2_x(\mathbb{R}^n)\to\s^2_x(\mathbb{R}^n) }\nonumber \\
\lesssim & \sum\limits_{j=0}^\infty \int_0^{2^{j/10}} dt 2^{(\delta-1) j} \times 1\times 2^{j/2 }t\nonumber\\
\lesssim & \sum\limits_{j=0}^\infty 2^{(\delta-3/10) j}\lesssim 1.
\end{align}
For $B_{52}$, we can use \eqref{1eq6} to get
\begin{align}
B_{52}\lesssim & \sum\limits_{j=0}^\infty \int_0^\infty dt 2^{j\delta}\|P^\pm e^{\pm itH_0}|P|^{1/2}\langle x\rangle^{-\sigma}\|_{\s^2_x(\mathbb{R}^n)\to \s^2_x(\mathbb{R}^n)}\nonumber\\
&\times\| \langle x\rangle^\sigma \frac{1}{|P|^{1/2}}F(|P|>2^{j/2})\langle x\rangle^{-\sigma}\|_{\s^2_x(\mathbb{R}^n)\to\s^2_x(\mathbb{R}^n)}\nonumber\\
\lesssim & \sum\limits_{j=0}^\infty 2^{(\delta-1/4)j}\lesssim 1.
\end{align}
So one has
\eq
B_5\lesssim 1.\label{B5}
\eeq
Based on \eqref{B3}, \eqref{B4} and \eqref{B5}, we get \eqref{1eq1}. We finish the proof.\\
\textbf{Proof of \eqref{100eq1}: } Let
\eq
B_3:=\int_0^\infty dt\|F(|x|\leq 1)\langle x\rangle^\delta P^\pm e^{\pm itH_0}\|_{\s^2_{x,\sigma}(\mathbb{R}^n)\cap \s^1_x(\mathbb{R}^n)\to\s^2_x(\mathbb{R}^n) },
\eeq
\eq
B_4:=\sum\limits_{j=0}^\infty\int_0^\infty \|F_{2^j}(|x|)\langle x\rangle^\delta P^\pm e^{\pm itH_0}\|_{\s^2_{x,\sigma}(\mathbb{R}^n)\cap \s^1_x(\mathbb{R}^n)\to\s^2_x(\mathbb{R}^n) } \chi(t\geq 2^{j/10})
\eeq
and
\eq
\tilde{B}_5:=\sum\limits_{j=0}^\infty\int_0^\infty \|F_{2^j}(|x|)\langle x\rangle^\delta P^\pm e^{\pm itH_0}\|_{\s^2_{x,\sigma}(\mathbb{R}^n)\cap H^{3/2}_x(\mathbb{R}^n)\to\s^2_x(\mathbb{R}^n) } \chi(t<2^{j/10}).
\eeq
Then
\eq
\text{(LHS) of }\eqref{100eq1}\leq B_3+B_4+\tilde{B}_5.
\eeq
If $\sigma>n/2$, then for sufficiently small $\epsilon>0$, we can see that
\eq
\sigma(1/2+\epsilon)-1>n/4-n\epsilon/2-1.
\eeq
Additionally, when $n\geq 45$, we can observe that
\eq
n/40-n\epsilon/20-1/10>0
\eeq
if $\epsilon>0$ is sufficiently small. So based on \eqref{B3}, \eqref{B4}, \eqref{B41} and \eqref{B42}, one has that when $\delta\in (0,n/40-n/20\times \epsilon-1/10)$, for sufficiently small $\epsilon>0$,
\eq
B_3\lesssim 1\label{lastB3}
\eeq
and
\eq
B_4\lesssim_{\delta,\sigma,\epsilon} 1.\label{lastB4}
\eeq
For $\tilde{B}_5$, break it into two pieces
\begin{align}
\tilde{B}_5\leq& \sum\limits_{j=0}^\infty \int_0^\infty dt \| \langle x\rangle^{\delta} F_{2^j}(|x|)P^\pm e^{\pm i tH_0}F(|P|\leq 2^{3j/4})\|_{\s^2_{x,\sigma}(\mathbb{R}^n)\cap H^{3/2}_x(\mathbb{R}^n)\to\s^2_x(\mathbb{R}^n) } \chi(t<2^{j/10})\nonumber\\
&+\sum\limits_{j=0}^\infty \int_0^\infty dt \| \langle x\rangle^\delta F_{2^j}(|x|)P^\pm e^{\pm i tH_0}F(|P|> 2^{3j/4})\|_{\s^2_{x,\sigma}(\mathbb{R}^n)\cap H^{3/2}_x(\mathbb{R}^n)\to\s^2_x(\mathbb{R}^n) } \chi(t<2^{j/10})\nonumber\\
=:&\tilde{B}_{51}+\tilde{B}_{52}.
\end{align}
Using Lemma \ref{PPpm}, we can see that for $\tilde{B}{51}$, if $N\in (22/3,\sigma)$ (where $\sigma>n/2\geq 45/2>22/3$) and $ \delta \in (1, 3N/20-1/10)$, then
\begin{align}
\tilde{B}_{51}\lesssim & \sum\limits_{j=0}^\infty \int_0^\infty dt 2^{(\delta-N) j} \| \langle x\rangle^N P^\pm \langle x\rangle^{-N}\|_{\s^2_x(\mathbb{R}^n)\to\s^2_x(\mathbb{R}^n) }\nonumber\\
&\times\| \langle x\rangle^N e^{\pm i tH_0}F(|P|\leq 2^{3j/4})\langle x\rangle^{-\sigma} \|_{\s^2_x(\mathbb{R}^n)\to\s^2_x(\mathbb{R}^n) } \chi(t<2^{j/10})\nonumber\\
\lesssim_N &\sum\limits_{j=0}^\infty 2^{(\delta-3N/20+1/10)j}\nonumber\\
\lesssim_{N,\delta}& 1.\label{B51}
\end{align}
For $\tilde{B}_{52}$, if $\delta\in (1,41/40)$, then we have
\begin{align}
\tilde{B}_{52}\leq & \sum\limits_{j=0}^\infty \int_0^\infty dt \| \langle x\rangle^\delta F_{2^j}(|x|)P^\pm e^{\pm itH_0}F(|P|> 2^{3j/4})\|_{ H^{3/2}_x(\mathbb{R}^n)\to\s^2_x(\mathbb{R}^n) } \chi(t<2^{j/10})\nonumber\\
\lesssim &\sum\limits_{j=0}^\infty \int_0^\infty dt \| \langle x\rangle^\delta F_{2^j}(|x|)P^\pm e^{\pm i tH_0}F(|P|> 2^{3j/4})\frac{1}{\langle P\rangle^{3/2}}\|_{ \s^2_x(\mathbb{R}^n)\to\s^2_x(\mathbb{R}^n) } \chi(t<2^{j/10})\nonumber\\
\lesssim& \sum\limits_{j=0}^\infty \frac{1}{2^{j(41/40-\delta)}}\nonumber\\
\lesssim_\delta & 1.\label{B52}
\end{align}
Using \eqref{B51}, \eqref{B52}, \eqref{lastB3}, and \eqref{lastB4}, we can see that if $\delta \in (1, \min\{n/40-0-1/10, 41/40, 3\sigma/20-1/10\})$, then
\eq
\text{(LHS) of }\eqref{100eq1}\leq B_3+B_4+\tilde{B}_5\lesssim_{\delta,\sigma} 1.
\eeq
\textbf{Proof of \eqref{Jan23.1}: } It follows from \eqref{1eq2}.
\end{proof}
Let
\begin{align}
C_\pm(t,T):=&\pm i\int_0^T ds P^{\pm}e^{\pm isH_0}F(|x|\geq10)V(x,t\pm s)U(t\pm s,t)\nonumber\\
&\mp \int_0^T ds P^{\pm}e^{\pm isH_0}[H_0,F(|x|\geq10)]U(t\pm s,t).
\end{align}
\begin{lemma}\label{Lem1}If $ \chi(|x|\geq 1)V(x,t)\psi(t)\in \s^\infty_t\s^2_{x,\sigma}(\mathbb{R}^{n+1})\cap \s^\infty_t\s^1_x(\mathbb{R}^{n+1})$ for some $\sigma> 2$ and if $\psi(t)$ satisfies \eqref{H1}, then for all $\delta\in(0,   \min\{ \sigma/20-1/10,1/40\})$, when $n\geq 5$,
\eq
\sup\limits_{t\geq 0}\sup\limits_{T\in[ 0,\infty]}\|\langle x\rangle^\delta C_\pm(t,T)\psi(t)\|_{\s^2_x(\mathbb{R}^n)}\lesssim_{\delta,E} 1,\quad j=1,2,3,4.\label{2eq1}
\eeq
\end{lemma}
\begin{proof}Based on the assumptions on $V$ and $\psi$, we have $  \chi(|x|\geq 1)V(x,t)\psi(t)\in \s^\infty_t\s^2_{x,\sigma}(\mathbb{R}^{n+1})\cap \s^\infty_t\s^1_x(\mathbb{R}^{n+1})$, $\chi(|x|\geq 1)V(x,t)e^{-itH_0}\psi(0)\in \s^\infty_t\s^2_{x,\sigma}(\mathbb{R}^{n+1})\cap \s^\infty_t\s^1_x(\mathbb{R}^{n+1})$, and in addition, we have
\eq
[H_0,F(|x|\geq 10)]\psi(t)=(H_0F)\times \psi-2\sum\limits_{j=1}^{n} \partial_{x_j}[F]\partial_{x_j}\psi(t)\in \s^\infty_t\s^2_{x,\sigma}(\mathbb{R}^{n+1})\cap \s^\infty_t\s^1_x(\mathbb{R}^{n+1}),
\eeq
and
\begin{align}
[H_0,F(|x|\geq 10)]\psi(t)=&(H_0F)\times \psi\nonumber\\
&-2\sum\limits_{j=1}^{n} \partial_{x_j}[F]\partial_{x_j}e^{-itH_0}\psi(0)\in \s^\infty_t\s^2_{x,\sigma}(\mathbb{R}^{n+1})\cap \s^\infty_t\s^1_x(\mathbb{R}^{n+1}).
\end{align}
Therefore, \eqref{2eq1} follows by using \eqref{1eq1}.
\end{proof}
Recall that
\eq
\psi_D(t)=\psi(t)-e^{-itH_0}\psi(0)
\eeq
and
\eq
V_D(x,t)=\mathcal{N}(|\psi_D(t)|, |x|,t ).
\eeq

 \begin{proof}[Proof of Lemma \ref{LempsiD}]\label{pfLempsiD}Using Duhamel's formula, one has
 \eq
 \psi_D(t)=-i\int_0^t ds e^{-i(t-s)H_0}V(x,s)\psi(s).
 \eeq
 The result follows from
 \begin{multline}
 \|\chi(t\geq 0)\langle x\rangle^{-2}\langle P\rangle^{3/2}\psi_D(t) \|_{\s^\infty_t\s^2_x(\mathbb{R}^{n+1})}\lesssim E+ \sum\limits_{j=1}^n\|\langle x\rangle^{-2}\langle P_j\rangle^{3/2}F(|P_j|>1)\psi_D(t) \|_{\s^\infty_t\s^2_x(\mathbb{R}^{n+1})}\\
 \lesssim E+ \sum\limits_{j=1}^n\int_0^{t}ds\|\langle x\rangle^{-2}\frac{\langle P_j\rangle^{3/2}F(|P_j|>1)}{P_j} e^{-i(t-s)H_0} [P_jV(x,s)\psi(s)] \|_{\s^\infty_t\s^2_x(\mathbb{R}^{n+1})}\\
 (\text{use Assumption \ref{asp4}})\lesssim_E 1
 \end{multline}
 and
 \begin{multline}
 \|\chi(t< 0)\langle x\rangle^{-2}\langle P\rangle^{3/2}\psi_D(t) \|_{\s^\infty_t\s^2_x(\mathbb{R}^{n+1})}\lesssim E+ \sum\limits_{j=1}^n\|\langle x\rangle^{-2}\langle P_j\rangle^{3/2}F(|P_j|>1)\psi_D(t) \|_{\s^\infty_t\s^2_x(\mathbb{R}^{n+1})}\\
 \lesssim E+ \sum\limits_{j=1}^n\int_0^{-t}ds\|\langle x\rangle^{-2}\frac{\langle P_j\rangle^{3/2}F(|P_j|>1)}{P_j} e^{-i(t+s)H_0} [P_jV(x,-s)\psi(-s)] \|_{\s^\infty_t\s^2_x(\mathbb{R}^{n+1})}\\
 (\text{use Assumption \ref{asp4}})\lesssim_E 1.
 \end{multline}
 Here we also use
 \eq
 \| \chi(|t\pm s|>1)\langle x\rangle^{-2} \frac{\langle P_j\rangle^{3/2}F(|P_j|>1)}{P_j} e^{-i(t\pm s)H_0}\langle x_j\rangle^{-2}\|_{\s^2_x(\mathbb{R}^n)\to\s^2_x(\mathbb{R}^n) }\lesssim \frac{\chi(|t\pm s|>1)}{\langle t\pm s \rangle^{2}}
 \eeq
 and
 \eq
 \int_0^{\mp t} \chi(|t\pm s|\leq 1)ds \| \langle x\rangle^{-2} e^{-i(t\pm s)H_0} \frac{\langle P_j\rangle^{3/2}F(|P_j|>1)}{P_j} \langle x_j\rangle^{-2}\|_{\s^2_x(\mathbb{R}^n)\to\s^2_x(\mathbb{R}^n) }\lesssim 1.
 \eeq
 We finish the proof.
 \end{proof}}
\section{Proof of Theorem \ref{thm} and examples}
\begin{proof}[Proof of Theorem \ref{thm}] Using incoming/outgoing decomposition, $F(|x|\geq 10)\psi(t)$ can be rewritten as
\eq
F(|x|\geq 10)\psi(t)=P^+F(|x|\geq 10)\psi(t)+P^-F(|x|\geq 10)\psi(t).
\eeq
Let
\eq
\psi_\pm :=\Omega_\pm^*\psi_0.
\eeq
Recall that
\begin{align}
C_\pm (t)=&P^\pm F(|x|\geq 10)-P^\pm \Omega_{t,\pm}^*\nonumber\\
=&\pm i\int_0^\infty ds P^\pm e^{\pm isH_0} F(|x|\geq 10)V(x,t\pm s)U(t\pm s,t)\nonumber\\
&+(\mp i)\int_0^\infty ds P^\pm e^{\pm is H_0} [H_0, F(|x|\geq 10)] U(t\pm s,t).
\end{align}
By approximating $P^\pm\psi(t)$ with $P^\pm e^{-itH_0}\psi_\pm$, and taking into account the intertwining property \eqref{inter}, we obtain
\begin{align}
F(|x|\geq 10)\psi(t)=&P^+\Omega_{t,+}^*\psi(t)+P^-\Omega_{t,+}^*\psi(t)+(C_+(t)+C_-(t))\psi(t)\nonumber\\
=&P^+e^{-itH_0}\psi_++P^-e^{-itH_0}\psi_-+(C_+(t)+C_-(t))\psi(t).
\end{align}
With the use of $P^\pm$, we can regard $\Omega_{t,\pm}^*\psi(t)$ as
\eq
\Omega_{t,\pm}^*\psi(t)=w\text{-}\lim\limits_{s\to \pm \infty} e^{isH_0}F(|x|\geq 10)\psi(t+s)\text{ on }\s^2_x(\mathbb{R}^n),
\eeq
which is demonstrated in Lemma \ref{Lem2}. See also Lemma \ref{Lem3}. Define
\eq
C(t):=C_+(t)+C_-(t).
\eeq
Then, let
\eq
\psi_{loc}(t):=C(t)\psi(t)+F(|x|<10)\psi(t).\label{oldloc}
\eeq
By applying Lemma \ref{Lem1}, we obtain that for all $\delta\in(0,   \min\{ \sigma/20-1/10,1/40\})$, when $n\geq 5$,
\eq
\| \langle x\rangle^\delta \psi_{loc}(t)\|_{\s^2_x(\mathbb{R}^n)}\lesssim_{\delta,E} \|\psi(0)\|_{\s^2_x(\mathbb{R}^n)}.
\eeq
Then we have
\begin{align}
\| \psi(t)-e^{-itH_0}\psi_+-\psi_{loc}(t)\|_{\s^2_x(\mathbb{R}^n)}=&\| P^+e^{-itH_0}\psi_++P^-e^{-itH_0}\psi_--e^{-itH_0}\psi_+\|_{\s^2_x(\mathbb{R}^n)}\nonumber\\
&+\| -P^-e^{-itH_0}\psi_++P^-e^{-itH_0}\psi_-\|_{\s^2_x(\mathbb{R}^n)}\to 0
\end{align}
as $t\to \infty$. This concludes the proof of Theorem \ref{thm}.

\end{proof}

\begin{proof}[Proof of Theorem \ref{thm2}]In the proof of Theorem \ref{thm}, we defined $\psi_{loc}(t)$ by \eqref{oldloc}. However, we now make some modifications:
\begin{align}
\psi_{loc}(t):=&\tilde{C}_+(t)\psi(0)+\tilde{C}_-(t)\psi(0)\nonumber\\
=:&\psi_{loc,+}(t)+\psi_{loc,-}(t)
\end{align}
where
\eq
V_{D}(x,t):=\mathcal{N}(|\psi_{D}(t)|,|x|,t)
\eeq
and
\begin{align}
\tilde{C}_\pm (t)\psi(0)=&\pm i\int_0^\infty ds P^\pm e^{\pm isH_0} V_{D}(x,t\pm s)\psi_D(x,t\pm s).
\end{align}
{\bf Localization property of $\psi_{loc}(t)$: }Due to \eqref{H1}, \eqref{con1}, Assumption \ref{asp2} and \eqref{1eq1} in Lemma \ref{out/in1}, we have $V_D(x,t)\in \s^\infty_t\s^2_{x,2}(\mathbb{R}^{n+1})$, which implies that
\eq
V_D(x,t)\psi_D(t)\in \s^\infty_t\s^1_{x,2}(\mathbb{R}^{n+1})\cap \s^\infty_t\s^2_{x,\sigma}(\mathbb{R}^{n+1})\text{ for some }\sigma>2
\eeq
and for all $\delta\in(0,   \min\{ \sigma/20-1/10,1/40\}) $, $n\geq 5$,
\eq
\|\langle x\rangle^\delta \tilde{C}_{\pm,1}(t)\psi(t)\|_{\s^2_x(\mathbb{R}^n)}\lesssim_{\delta,E} 1.
\eeq
So for all $\delta\in(0,   \min\{ \sigma/20-1/10,1/40\}) $,
\eq
\|\langle x\rangle^\delta \psi_{loc}(t)\psi(t)\|_{\s^2_x(\mathbb{R}^n)}\lesssim_{\delta,E} 1.
\eeq
{\bf Asymptotic decomposition:} This modification will not change the validity of \eqref{eq3}:
\begin{enumerate}
\item Let
\eq
\begin{split}
C_{\pm,1}(t)\psi(t):=&\pm i\int_0^\infty ds P^\pm e^{\pm isH_0}V(x,t\pm s)\psi(t\pm s).
\end{split}
\eeq
Due to \eqref{con1}, \eqref{1eq1} in Lemma \ref{out/in1} and the fact that $V(x,t)\in \s^\infty_t\s^2_{x,2}(\mathbb{R}^{n+1})$, we have
\eq
V(x,t)\psi(t)\in \s^\infty_t\s^1_{x,2}(\mathbb{R}^{n+1})\cap \s^\infty_t\s^2_{x,\sigma}(\mathbb{R}^{n+1})\text{ for some }\sigma>2
\eeq
and for all $\delta\in(0,   \min\{ \sigma/20-1/10,1/40\}) $,
\eq
\|\langle x\rangle^\delta C_{\pm,1}(t)\psi(t)\|_{\s^2_x(\mathbb{R}^n)}\lesssim_{\delta,E} 1.
\eeq
If we take
\eq
\tilde{\psi}_{loc}(t):=C_{+,1}(t)\psi(t)+C_{-,1}(t)\psi(t),
\eeq
we have
\begin{align}
\psi(t)=& P^+\psi(t)+P^-\psi(t)\nonumber\\
=&P^+\Omega_{t,+}^*\psi(t)+P^-\Omega_{t,-}^*\psi(t)+C_{+,1}(t)\psi(t)+C_{-,1}(t)\psi(t)\nonumber\\
=&P^+e^{-itH_0}\Omega_+^*\psi(0)+P^-e^{-itH_0}\Omega_-^*\psi(0)+\tilde{\psi}_{loc}(t)
\end{align}
which implies
\begin{align}
& \| \psi(t)-e^{-itH_0}\Omega_+^*\psi(0)-\tilde{\psi}_{loc}(t)\|_{\s^2_x(\mathbb{R}^n)}\nonumber\\
=&\| -P^-e^{-itH_0}\Omega_+^*\psi(0)+P^-e^{-itH_0}\Omega_-^*\psi(0)\|_{\s^2_x(\mathbb{R}^n)}\nonumber\\
\to & 0\label{Jan22.1}
\end{align}
as $t\to \infty$.
\item
Due to \eqref{con2} and \eqref{1eq3} in Lemma \ref{out/in1}, because for any $\sigma'\in(0,1/2)$,
\eq
\| \psi(t)-\psi_D(t)\|_{\s^2_{x,-\sigma'}}=\| e^{-itH_0}\psi(0)\|_{\s^2_{x,-\sigma'}}\to 0
\eeq
as $t\to \infty$, we have
\begin{align}
\|C_{\pm,1}(t)\psi(t)-\tilde{C}_\pm (t)\psi(0)\|_{\s^2_x(\mathbb{R}^n)}=&\nonumber\\
\|\pm i\int_0^\infty ds P^\pm e^{\pm isH_0} (V(x,t\pm s)&\psi(t\pm s)-V_D(x,t\pm s)\psi_D(t\pm s))\|_{\s^2_x(\mathbb{R}^n)}\nonumber\\
\leq \|\pm i\int_0^{t/2} ds P^\pm e^{\pm isH_0} (V(x,t\pm s)&\psi(t\pm s)-V_D(x,t\pm s)\psi_D(t\pm s))\|_{\s^2_x(\mathbb{R}^n)}\nonumber\\
+\|\pm i\int_{t/2}^\infty ds P^\pm e^{\pm isH_0} (V(x,t\pm s)&\psi(t\pm s)-V_D(x,t\pm s)\psi_D(t\pm s))\|_{\s^2_x(\mathbb{R}^n)}\nonumber\\
\lesssim_{E} \sup\limits_{s\geq t/2} \| e^{-isH_0}\psi(0)\|_{\s^2_{x,-\sigma'}} +\int_{t/2}^\infty ds& \|P^\pm e^{\pm i sH_0}\|_{\s^2_{x,\sigma}(\mathbb{R}^n)\cap \s^1_x(\mathbb{R}^n)\to \s^2_x(\mathbb{R}^n)}\to 0\label{Jan22.2}
\end{align}
as $t\to \infty$. According to \eqref{Jan22.1} and \eqref{Jan22.2}, by setting
\eq
\psi_{loc}(t)=\tilde{C}_+ (t)\psi(0)+\tilde{C}_- (t)\psi(0),
\eeq
we have
\begin{align}
& \| \psi(t)-e^{-itH_0}\Omega_+^*\psi(0)-\psi_{loc}(t)\|_{\s^2_x(\mathbb{R}^n)}\nonumber\\
\leq &\| \psi(t)-e^{-itH_0}\Omega_+^*\psi(0)-\tilde{\psi}_{loc}(t)\|_{\s^2_x(\mathbb{R}^n)}\nonumber\\
&+\|C_{+,1}(t)\psi(t)-\tilde{C}_+ (t)\psi(0)\|_{\s^2_x(\mathbb{R}^n)}+\|C_{-,1}(t)\psi(t)-\tilde{C}_- (t)\psi(0)\|_{\s^2_x(\mathbb{R}^n)} \nonumber\\
\to & 0
\end{align}
as $t\to \infty$.
\end{enumerate}
Let
\begin{align*}
\psi_{loc}(t)=&\tilde{C}_+ (t)\psi(0)+\tilde{C}_- (t)\psi(0)\\
=:&\psi_{loc,1}(t)+\psi_{loc,2}(t).
\end{align*}
{\bf Estimate for $A^2\psi_{loc,1}(t)$:} For $A^2\psi_{loc,1}(t)$, we have
\begin{subequations}
    \eq
    i[H_0,A]=2H_0;
    \eeq
    \eq
    (-i)[H_0,A^2]=-4H_0A-4iH_0;
    \eeq
    \eq
    (-i)[H_0, (-i)[H_0,A^2]]=8H_0;
    \eeq
    \begin{align}
    A^2e^{isH_0}=&e^{isH_0}A^2+(-i)\int_0^s du e^{-iuH_0}[H_0,A^2]e^{iuH_0}\nonumber\\
    =&e^{isH_0}A^2+(-i)s[H_0,A^2]+(-i)\int_0^sdu\int_0^udv e^{-ivH_0}(-i)[H_0,[H_0,A^2]]e^{ivH_0}\nonumber\\
    =&e^{isH_0}A^2+(-4H_0A-4iH_0)s+4H_0^2s^2,
    \end{align}
\end{subequations}
which implies that
\begin{align}
A^2\psi_{loc,1}(t)=&i\int_0^\infty ds P^+ A^2 e^{ i sH_0}V_D(x,t+ s)\psi_D(t+ s)\nonumber\\
=&i\int_0^\infty ds P^+\left(4s^2e^{isH_0}(-\Delta)^2-4se^{isH_0}(-\Delta)A+(-4is)e^{isH_0}(-\Delta)A\right.\nonumber\\
&\left.+e^{isH_0}A^2\right)\times V_D(x,t+ s)\psi_D(t+ s).
\end{align}
Break $A^2\psi_{loc,1}(t)$ into two pieces
\begin{align}
A^2\psi_{loc,1}(t)=&i\int_0^\infty ds P^+F(|P|\leq 1)\left(4s^2e^{isH_0}(-\Delta)^2-4se^{isH_0}(-\Delta)A+\right.\nonumber\\
&\left.(-4is)e^{isH_0}(-\Delta)A+e^{isH_0}A^2\right)V_D(x,t+ s)\psi_D(t+ s)\nonumber\\
&+i\int_0^\infty ds P^+F(|P|> 1)\left(4s^2e^{isH_0}(-\Delta)^2-4se^{isH_0}(-\Delta)A\right.\nonumber\\
&\left.+(-4is)e^{isH_0}(-\Delta)A+e^{isH_0}A^2\right)V_D(x,t+ s)\psi_D(t+ s)\nonumber\\
=:&C_{m,+,1}(t)\psi(0)+C_{m,+,2}(t)\psi(0).
\end{align}
For $C_{m,+,1}(t)\psi(0)$, according to \eqref{Jan23.1} in Lemma \ref{out/in1}, one has that when $\sigma>6$, $\sigma+a-2>4$ and $V_D(x,t)\in \s^2_t\s^2_{x,2}(\R^{n+1})$, for all $a=0,1,2,$, we have
\begin{align}
\|C_{m,+,1}(t)\psi(0)\|_{\s^2_x(\mathbb{R}^n)}\lesssim &\sum\limits_{a=0}^2 \int_0^\infty ds s^a\| P^+e^{isH_0}(-\Delta)^aF(|P|\leq 1)\|_{\s^2_{x,\sigma+a-2}(\mathbb{R}^n)\cap \s^1_x(\mathbb{R}^n)\to\s^2_x(\mathbb{R}^n) }\nonumber\\
&\times\|  F(|P|\leq 10) A^{2-a} V_D(x,t+s)\psi_D(t+s)\|_{\s^2_{x,\sigma+a-2}(\mathbb{R}^n)\cap \s^1_x(\mathbb{R}^n)}\nonumber\\
\lesssim & \|V_D(x,t+s)\psi_D(t+s)\|_{\s^2_{x,\sigma}(\mathbb{R}^n)\cap \s^1_{x,2}(\mathbb{R}^n)} \nonumber\\
\lesssim_{E,\|V(x,t)\|_{\s^\infty_t\s^2_{x,2}(\R^{n})}}& 1.\label{Cm+1}
\end{align}
Here we also use that \eqref{H1} implies
\eq
\sup\limits_{t\in \mathbb{R}}\| \psi_D(t)\|_{H^1_x(\mathbb{R}^n)}\leq 2E<\infty
\eeq
and therefore
\eq
\| V_D(x,t+s)\psi_D(t+s)\|_{\s^2_{x,\sigma}(\mathbb{R}^n)\cap \s^1_{x,2}(\mathbb{R}^n)}\lesssim_{E,\|V(x,t)\|_{\s^\infty_t\s^2_{x,2}(\R^{n})}}1
\eeq
by using \eqref{con1}. For $C_{m.+,2}(t)\psi(0)$, according to \eqref{Sep29.1} and \eqref{Oct.1} in Lemma \ref{out/in1}, using Lemma \ref{LempsiD} and \eqref{con5}, we have that for some $\epsilon>0$ close to $0$,
\begin{align}
\|C_{m,+,2}(t)\psi(0)\|_{\s^2_x(\mathbb{R}^n)}\lesssim &\sum\limits_{a=0}^2 \int_0^\infty ds s^{a}\| P^+e^{isH_0}F(|P|>1)|P|^{a+1-\epsilon}\|_{\s^2_{x,3+\epsilon}(\mathbb{R}^n)\to\s^2_x(\mathbb{R}^n) }\nonumber\\
&\times\| |P|^{a-1+\epsilon} A^{2-a}V_D(x,t+s)\psi_D(t+s)\|_{\s^2_{x,3+\epsilon}(\mathbb{R}^n)}\nonumber\\
\lesssim &\sum\limits_{a=0}^2 \int_0^\infty ds s^{a}\| P^+e^{isH_0}|P|^{a+1-\epsilon}\|_{\s^2_{x,3+\epsilon}(\mathbb{R}^n)\to\s^2_x(\mathbb{R}^n) }\nonumber\\
&\times \| \langle P\rangle^{3/2}\langle x\rangle^{5+\epsilon} V_D(x,t+s)\psi_D(t+s)\|_{\s^2_{x}(\mathbb{R}^n)} \nonumber\\
\lesssim_{\epsilon,E} &1,\label{Cm+2}
\end{align}
where we use
$$
\| \langle x\rangle^{3+\epsilon} |P|^{a-1+\epsilon} A^{2-a}\langle x\rangle^{-(5+\epsilon)} \langle P\rangle^{-3/2}\|_{\s^2_{x}(\mathbb{R}^n)\to\s^2_{x}(\mathbb{R}^n)}\lesssim 1.
$$
Based on \eqref{Cm+1} and \eqref{Cm+2}, one has
\eq
\|A^2\psi_{loc,1}(t)\|_{\s^2_x(\mathbb{R}^n)}\lesssim_E 1.\label{Aloc1}
\eeq
Similarly, one has
\eq
\|A^2\psi_{loc,2}(t)\|_{\s^2_x(\mathbb{R}^n)}\lesssim_E 1.\label{Aloc2}
\eeq
Based on \eqref{Aloc1} and \eqref{Aloc2}, one has
\eq
\| A^2 \psi_{loc}(t)\|_{\s^2_x(\mathbb{R}^n)}\lesssim_E 1.
\eeq
We finish the proof of \eqref{smooth}. When $n\geq 45$ and $\sigma>n/2$, due to Assumption \ref{asp3},
\eq
\| \langle P \rangle^{3/2} V_D(x,t)\psi_D(t)\|_{\s^2_x(\mathbb{R}^n)}\lesssim \| \langle P \rangle^{3/2} \langle x\rangle^6V_D(x,t)\psi_D(t)\|_{\s^2_x(\mathbb{R}^n)}\lesssim_E1.
\eeq
So \eqref{locallast} follows by using \eqref{100eq1} in Lemma \ref{out/in1}.

\end{proof}
\begin{example}[Example of Theorem \ref{thm}]\label{example}When $n=5$ and $\mathcal{N}=\pm \lambda |\psi|^p$ for any $p\in (1,\frac{4}{3}]$, \eqref{condition} is satisfied and we have \eqref{eq2}.
\end{example}
\begin{proof}Using \eqref{H1}, we can conclude that
\eq
\sup\limits_{t\in \mathbb{R}}\|\psi(t)\|_{\s^{10/3}_x(\mathbb{R}^5)}\lesssim \sup\limits_{t\in \mathbb{R}}\|\psi(t)\|_{H^1_x(\mathbb{R}^5)}<\infty.
\eeq
Therefore,
\eq
\| |\psi(t)|^{p+1}\|_{\s^\infty_t\s^1_x(\mathbb{R}^{5+1})}<\infty
\eeq
and
\eq
\| |\psi(t)|^{p}\|_{\s^\infty_t\s^2_x(\mathbb{R}^{5+1})}<\infty
\eeq
for all $p\in (1,4/3]$. Furthermore, we have that for all $p\in (1,4/3]$,
\eq
\|\chi(|x|\geq 1) |\psi(t)|^{p+1}\|_{\s^2_{x,\sigma}(\mathbb{R}^{5+1})}<\infty
\eeq
for some $\sigma>2$, since \eqref{H1} implies
\eq
|\psi(t) |\lesssim_E \frac{1}{|x|^2}\text{ when }|x|\geq 1.
\eeq
This completes the proof.
\end{proof}
\begin{example}[Example of Theorem \ref{thm}]\label{example0}When $n=5$, $\mathcal{N}=W(x,t), W(x,t)\pm \lambda |\psi|^p$ for any $p\in (1,\frac{4}{3}]$ and some $W$ satisfying that for some $\sigma>2$,
\eq
\begin{cases}
W(x,t)\in \s^\infty_t\s^2_x(\mathbb{R}^{5+1})\\
\chi(|x|\geq 1)W(x,t)\in \s^\infty_t\s^{\infty}_{x,\sigma}(\mathbb{R}^{5+1})\cap\s^\infty_t\s^{2}_x(\mathbb{R}^{5+1})
\end{cases}.
\eeq
\end{example}
\begin{proof}When $\mathcal{N}=W(x,t)$, it follows directly because $W$ is a linear interaction. When $\mathcal{N}=W(x,t)\pm \lambda |\psi|^p $, it follows from the proof of Example \ref{example}.
\end{proof}
\begin{example}[Example of Theorem \ref{thm2}]\label{example2}When $n=45$, $\mathcal{N}=- \lambda |\frac{\psi|^p}{1+|\psi|^p}$ for $p>3, \lambda>0$, we have \eqref{eq2} and \eqref{smooth}.
\end{example}
\begin{proof} $\mathcal{N}(|\psi(t)|,|x|,t)\in \s^\infty_t\s^2_{x,2}(\mathbb{R}^{n+1})$: Using
\begin{align}
\frac{|\psi|^{p}}{1+|\psi|^p}\leq& 1,
\end{align}
and 
\eq
|\psi(x,t)|\lesssim_{n,E} \frac{1}{|x|^{(n-1)/2}},\quad |x|\geq1,
\eeq
we have
\eq
\begin{split}
\| \lambda \frac{|\psi|^{p}}{1+|\psi|^p} \|_{\s^\infty_t\s^2_{x,2}(\mathbb{R}^{n+1})}\leq &\| \chi(|x|\leq 1)\lambda \frac{|\psi|^{p}}{1+|\psi|^p} \|_{\s^\infty_t\s^2_{x,2}(\mathbb{R}^{n+1})}+\|\chi(|x|> 1) \lambda \frac{|\psi|^{p}}{1+|\psi|^p} \|_{\s^\infty_t\s^2_{x,2}(\mathbb{R}^{n+1})}\\
\lesssim_E& \lambda\left(1+\| \psi(t)\|_{\s^\infty_t\s^2_x(\mathbb{R}^{n+1})}\right)\\
\lesssim_E &\lambda \| \psi(0)\|_{\s^2_x(\mathbb{R}^n)}.
\end{split}
\eeq
Assumption \ref{asp1}: \eqref{H1} implies
\eq
|\psi(t) |\lesssim_E \frac{1}{|x|^{\frac{n-1}{2}}}\text{ when }|x|\geq 1.
\eeq
Therefore,
\eq
|\lambda \frac{|\psi|^{p}}{1+|\psi|^p} |\lesssim_E \lambda \times \frac{1}{\langle x\rangle^{(n-1)p/2}}\label{asp1eq1}
\eeq
with $(n-1)p/2>3(n-1)/2>n/2$. So Assumption \ref{asp1} is satisfied.\\
Assumption \ref{asp4}: Let
\eq
g(k):=\frac{d}{dk}[\frac{k^{p/2}}{1+k^{p/2}}]=\frac{p/2\times k^{p/2-1}}{1+k^{p/2}}-\frac{k^{p/2}\times p/2\times k^{p/2-1}}{(1+k^{p/2})^2},\quad k\geq 0.
\eeq
Therefore,
\eq
|g(k)|\lesssim_p \frac{|k|^{p/2-1}}{1+|k|^{p/2}}\quad\text{ for all }k\geq 0.
\eeq
Compute $\partial_{x_j}[ \frac{|\psi|^{p}}{1+|\psi|^p}\psi]$, $j=1,\cdots,n$:
\begin{align}
\partial_{x_j}[ \frac{|\psi|^{p}}{1+|\psi|^p}\psi]=&\frac{|\psi|^{p}}{1+|\psi|^p}\partial_{x_j}[\psi]+g(k)\vert_{k=|\psi|^2} (\partial_{x_j}[\psi^*] \psi+\psi^* \partial_{x_j}[\psi])\times \psi,\label{Pjpsi}
\end{align}
which implies that
\begin{align}
\| \partial_{x_j}[ \frac{|\psi|^{p}}{1+|\psi|^p}\psi]\|_{\s^\infty_t\s^2_{x,2}(\mathbb{R}^{n+1})}\lesssim_p&\| k^p\vert_{k=|\psi|}\times (|\partial_{x_j}[\psi]|+|\partial_{x_j}[\psi^*]|)\|_{\s^\infty_t\s^2_{x,2}(\mathbb{R}^{n+1})}\nonumber\\
\lesssim_{p,E}& \| \frac{1}{\langle x\rangle^{2p}}  (|\partial_{x_j}[\psi]|+|\partial_{x_j}[\psi^*]|)\|_{\s^\infty_t\s^2_{x,2}(\mathbb{R}^{n+1})}\nonumber\\
\lesssim_{p,E}& 1.
\end{align}
Therefore, Assumption \ref{asp4} is satisfied.\\
Assumption \ref{asp2}: Write $V(x,t)\psi(t)-V_D(x,t)\psi_D(t)$ as
\begin{align}
V(x,t)\psi(t)-V_D(x,t)\psi_D(t)=&V(x,t)e^{-itH_0}\psi(0)+(V(x,t)-V_D(x,t))\psi_D(t)\nonumber\\
=:&V\psi_1(t)+V\psi_2(t).
\end{align}
For $V\psi_1(t)$, one has that when $ \sigma'\in (0,\frac{n-1}{2}\times p-n/2-2)$ and $\sigma\in (\frac{n}{2},\frac{n-1}{2}\times p-\frac{n}{2}-\sigma')$,
\eq
\frac{1}{\langle x\rangle^{\frac{n-1}{2}\times p-(\sigma+\sigma')}}\in \s^\infty_x(\mathbb{R}^n), \quad \frac{1}{\langle x\rangle^{\frac{n-1}{2}\times p-\sigma'}}\in \s^2_{x,2}(\mathbb{R}^n),
\eeq
and therefore using H\"older's inequality,
\begin{align}
\| V\psi_1(t)\|_{\s^2_{x,\sigma}(\mathbb{R}^n)\cap \s^1_{x,2}(\mathbb{R}^n)}\leq &\| V(x,t)\|_{\s^\infty_{x,\sigma+\sigma'}\cap \s^2_{x,\sigma'+2}} \| e^{-itH_0}\psi(0)\|_{\s^2_{x,-\sigma'}}\nonumber\\
(\text{use \eqref{asp1eq1}})\lesssim_E \| e^{-itH_0}\psi(0)\|_{\s^2_{x,-\sigma'}}.
\end{align}
For $V\psi_2(t)$, let
\eq
\tilde{g}(k):=\frac{d}{dk}[\frac{-\lambda k^p}{1+k^p}],\quad k\geq0.
\eeq
Then
\eq
| \tilde{g}(k)|\lesssim_{p,\lambda} \frac{k^{p-1}}{1+k^p},\quad k\geq 0.
\eeq
Using
\eq
||f|-|g||\leq  |f-g|,
\eeq
based on the fundamental theorem of calculus, one has
\begin{align}
\chi(|x|\leq 1)| V\psi_2(t)|=&\chi(|x|\leq 1) |  \int_0^{|\psi(t)|-|\psi_D(t)|}dk \tilde{g}(|\psi_D(t)|+k) | |\psi_D(t)|\nonumber\\
\lesssim_{p,\lambda}  & \chi(|x|\leq 1)|  \int_0^{|\psi(t)|-|\psi_D(t)|}dk \frac{1}{1+|\psi_D|+k}    |\psi_D(t)|\nonumber\\
\lesssim_{p,\lambda}& \chi(|x|\leq 1) |e^{-itH_0}\psi(0)|^{2/n}|\psi_D(t)|,
\end{align}
which implies that by using H\"older's inequality, for any $\sigma',\sigma>0$,
\begin{align}
\| \chi(|x|\leq 1)V\psi_2(t)\|_{\s^2_{x,\sigma}(\mathbb{R}^n)}\lesssim_{p,\lambda}&\|e^{-itH_0}\psi(0) \|_{\s^2_{x,-\sigma'}}^{2/n}\|  \psi_D(t)\|_{\s^{2n/(n-2)}_x(\mathbb{R}^n)}\nonumber\\
\lesssim_{p,\lambda,E}& \|e^{-itH_0}\psi(0) \|_{\s^2_{x,-\sigma'}(\mathbb{R}^n)}^{2/n}.
\end{align}
We also have
\begin{align}
\chi(|x|\leq 1)| V\psi_2(t)|=&\chi(|x|\leq 1) |  \int_0^{|\psi(t)|-|\psi_D(t)|}dk \tilde{g}(|\psi_D(t)|+k) | |\psi_D(t)|\nonumber\\
\lesssim_{\lambda,p}  & \chi(|x|\leq 1)| |e^{-itH_0}\psi(0)|    |\psi_D(t)|,\label{sVpsi1}
\end{align}
which implies that by using H\"older's inequality, for any $\sigma'>0$,
\begin{align}
\| \chi(|x|\leq 1)V\psi_2(t)\|_{\s^1_{x,2}(\mathbb{R}^n)}\lesssim_{\lambda,p}&\|e^{-itH_0}\psi(0) \|_{\s^2_{x,-\sigma'}}\|  \psi_D(t)\|_{\s^2_x(\mathbb{R}^n)}\nonumber\\
\lesssim_{\lambda,p,E} &\|e^{-itH_0}\psi(0) \|_{\s^2_{x,-\sigma'}(\mathbb{R}^n)}.\label{sVpsi2}
\end{align}
Hence, based on \eqref{sVpsi1} and \eqref{sVpsi2}, one has that for any $\sigma'>0$,
\eq
\| \chi(|x|\leq 1)V\psi_2(t)\|_{\s^2_{x,\sigma}(\mathbb{R}^n)\cap\s^1_{x,2}(\mathbb{R}^n)}\lesssim_{p,\lambda,E} \|e^{-itH_0}\psi(0) \|_{\s^2_{x,-\sigma'}(\mathbb{R}^n)}^{2/n}+\|e^{-itH_0}\psi(0) \|_{\s^2_{x,-\sigma'}(\mathbb{R}^n)}.\label{sVpsi0}
\eeq
For $\chi(|x|>1)V\psi_2(t)$, we have
\begin{align}
\chi(|x|> 1)| V\psi_2(t)|=&\chi(|x|> 1) |  \int_0^{|\psi(t)|-|\psi_D(t)|}dk \tilde{g}(|\psi_D(t)|+k) | |\psi_D(t)|\nonumber\\
\lesssim_{\lambda,p}  & \chi(|x|>1)| |e^{-itH_0}\psi(0)| \times\left(|\psi(t)|^{p-1}+|\psi_D(t)|^{p-1} \right)   |\psi_D(t)|,\label{lVpsi1}
\end{align}
which implies that when $ \sigma'\in (0,\frac{n-1}{2}\times p-n/2-2)$ and $\sigma\in (\frac{n}{2},\frac{n-1}{2}\times p-\frac{n}{2}-\sigma')$,
\eq
\| \chi(|x|> 1)V\psi_2(t)\|_{\s^2_{x,\sigma}(\mathbb{R}^n)\cap \s^1_{x,2}(\mathbb{R}^n)}\lesssim_{\lambda,p,E}\|e^{-itH_0}\psi(0)\|_{\s^2_{x,-\sigma'}(\mathbb{R}^n)}.
\eeq
Thus, Assumption \ref{asp2} is satisfied.\par
Assumption \ref{asp3}: Using the chain rule, we have
\eq
\begin{split}
\| \langle P\rangle^{3/2}\langle x\rangle^6 V_D(x,t)\psi_D(t)\|_{\s^2_x(\mathbb{R}^n)}\lesssim_E&1+\sum\limits_{j=1}^n \| F(|P_j|>1)\langle P_j\rangle^{3/2}\langle x\rangle^6 V_D(x,t)\psi_D(t)\|_{\s^2_x(\mathbb{R}^n)}.
\end{split}
\eeq
Since
\eq
\begin{split}
\partial_{x_j}[ \frac{|\psi_D|^{p}}{1+|\psi_D|^p}\psi_D]=&\frac{|\psi_D|^{p}}{1+|\psi_D|^p}\partial_{x_j}[\psi_D]+g(k)\vert_{k=|\psi_D|^2} (\partial_{x_j}[\psi_D^*] \psi_D+\psi_D^* \partial_{x_j}[\psi_D])\times \psi_D,
\end{split}\label{PjpsiD}
\eeq
we have
\begin{multline}
\| F(|P_j|>1)\langle P_j\rangle^{3/2}\langle x\rangle^6 V_D(x,t)\psi_D(t)\|_{\s^2_x(\mathbb{R}^n)}\leq\\
\|F(|P_j|>1)\frac{\langle P_j\rangle^{3/2}}{iP_j} [\partial_{x_j}[\langle x\rangle^6]]V_D(x,t)\psi_{D}(x,t) \|_{\s^2_x(\mathbb{R}^n)}+\\
\|F(|P_j|>1)\frac{\langle P_j\rangle^{3/2}}{iP_j} f_{D,1}(x,t) \partial_{x_j}[\psi_{D}(x,t)] \|_{\s^2_x(\mathbb{R}^n)}+\\
\|F(|P_j|>1)\frac{\langle P_j\rangle^{3/2}}{iP_j} f_{D,2}(x,t) \partial_{x_j}[\psi^*_{D}(x,t)] \|_{\s^2_x(\mathbb{R}^n)}
\end{multline}
where
\eq
f_{D,1}(x,t):=\lambda\langle x\rangle^6 ( \frac{|\psi_D|^{p}}{1+|\psi_D|^p}+kg(k)\vert_{k=|\psi_D|^2})\in \s^\infty_tW^{1,\infty}_{x}(\mathbb{R}^{n+1})\cap \s^\infty_t\s^\infty_{x,2}(\mathbb{R}^{n+1})\label{fD1}
\eeq
and
\eq
f_{D,2}(x,t):=\lambda\langle x\rangle^6 g(k)\vert_{k=|\psi_D|^2}\psi_D^2\in \s^\infty_tW^{1,\infty}_{x}(\mathbb{R}^{n+1})\cap \s^\infty_t\s^\infty_{x,2}(\mathbb{R}^{n+1}).\label{fD2}
\eeq
Using the chain rule, Assumption \ref{asp4} and Assumption \ref{asp1}, we have
\begin{multline}
\|F(|P_j|>1)\frac{\langle P_j\rangle^{3/2}}{iP_j} [\partial_{x_j}[\langle x\rangle^6]]V_D(x,t)\psi_{D}(x,t) \|_{\s^2_x(\mathbb{R}^n)}\leq \\
\|  [\langle P_j\rangle^{1/2}\partial_{x_j}[\langle x\rangle^6]]V_D(x,t)\psi_{D}(x,t)  \|_{\s^2_x(\mathbb{R}^n)}+\|  \partial_{x_j}[\langle x\rangle^6] \langle P_j\rangle^{1/2}[V_D(x,t)\psi_{D}(x,t)]  \|_{\s^2_x(\mathbb{R}^n)}\\
\lesssim_E 1.\label{ex2:1}
\end{multline}
Using Lemma \ref{LempsiD}, the chain rule, Assumption \ref{asp4} and Assumption \ref{asp1} , we have
\begin{multline}
\|F(|P_j|>1)\frac{\langle P_j\rangle^{3/2}}{iP_j} f_{D,1}(x,t) \partial_{x_j}[\psi_{D}(x,t)] \|_{\s^2_x(\mathbb{R}^n)}\leq \\
\| \langle P_j\rangle^{1/2} f_{D,1}(x,t)\|_{\s^\infty_x(\mathbb{R}^5)}\|\partial_{x_j}[\psi_{D}(x,t)]\|_{\s^2_x(\mathbb{R}^n)}+\|  f_{D,1}(x,t)\|_{\s^\infty_{x,2}(\mathbb{R}^n)}\| \langle P_j\rangle^{1/2}\partial_{x_j}[\psi_{D}(x,t)]\|_{\s^2_{x,-2}(\mathbb{R}^n)}\\
(\text{use \eqref{fD1}})\lesssim_E 1\label{ex2:2}
\end{multline}
and
\begin{multline}
\|F(|P_j|>1)\frac{\langle P_j\rangle^{3/2}}{iP_j} f_{D,2}(x,t) \partial_{x_j}[\psi_{D}^*(x,t)] \|_{\s^2_x(\mathbb{R}^n)}\leq\\
 \| \langle P_j\rangle^{1/2} f_{D,2}(x,t)\|_{\s^\infty_x(\mathbb{R}^n)}\|\partial_{x_j}[\psi_{D}^*(x,t)]\|_{\s^2_x(\mathbb{R}^n)}+\|  f_{D,2}(x,t)\|_{\s^\infty_{x,2}(\mathbb{R}^n)}\| \langle P_j\rangle^{1/2}\partial_{x_j}[\psi_{D}^*(x,t)]\|_{\s^2_{x,-2}(\mathbb{R}^n)}\\
(\text{use \eqref{fD2}})\lesssim_E 1.\label{ex2:3}
\end{multline}
According to \eqref{ex2:1}, \eqref{ex2:2} and \eqref{ex2:3}, we have
\eq
\| F(|P_j|>1)\langle P_j\rangle^{3/2}\langle x\rangle^6 V_D(x,t)\psi_D(t)\|_{\s^2_x(\mathbb{R}^n)}\lesssim_E 1,
\eeq
which implies
\eq
\| \langle P\rangle^{3/2}\langle x\rangle^6 V_D(x,t)\psi_D(t)\|_{\s^2_x(\mathbb{R}^n)}\lesssim_E1.
\eeq
Therefore, Assumption \ref{asp3} is satisfied. We finish the proof.
\end{proof}
\section{Application}
 In this section, we present an application of this method, demonstrating that for certain nonlinear Schr\"odinger equations with a linear potential that is quasi-periodic in time, global Strichartz estimates hold in $t$ for scattering states. Consider the linear time-dependent Schr\"odinger equation:
\eq
\begin{cases}
    i\partial_t\psi(x,t)=(-\Delta_x+V_0(x)+V_{s}(x,t)+V_{c}(x,t))\psi(x,t)\\
    \psi(x,0)=\psi_0\in \s^2(\R^5)
\end{cases}, \quad (x,t)\in \R^{5+1}.\label{linear}
\eeq
Here, $V_0(x), V_s(x,t)$ and $V_c(x,t)$ are real valued and satisfy following assumption:
\begin{assumption}\label{dasp1}Assume $V_0(x)\in \s^\infty_{x,\sigma}(\R^5)$ for some $\sigma>6$, and $V_{s}(x,t)$ and $V_{c}(x,t)$ have the form: 
\begin{align}
V_{s}(x,t)=&\sum\limits_{j=1}^{N_1}\sin(\omega_{sj}t)V_{sj}(x), 
\end{align}
and
\begin{align}
V_{c}(x,t)=&\sum\limits_{j=1}^{N_2}\cos(\omega_{cj}t)V_{cj}(x), 
\end{align}
where $V_{sj}(x),V_{cj}(x)\in \s^\infty_\sigma(\R^5)$ for $\sigma>6$, and $\omega_{sj}$ are pairwise irrational, as are $\omega_{cj}$. 
    
\end{assumption}
According to \cite{SW20221}, we have 
\eq
\Omega_\alpha^*:=s\text{-}\lim\limits_{t\to \infty} e^{itH_0}F_c(\frac{|x-2tP|}{t^\alpha}\leq 1)U_{lin}(t,0)\text{ exists on }\s^2_x(\R^5)
\eeq
for $\alpha\in (0,3/5)$, and 
\eq
P_{sc}(0):=s\text{-}\lim\limits_{t\to \infty}U_{V}(0,t)F_c(|x-2tP|\leq t^\alpha)U_{V}(t,0),\text{ exists on }\s^2_x(\R^5),
\eeq
where $U_{V}(t,0)$ denotes the solution operator to system \eqref{linear} and $F_c$ denotes a smooth characteristic function. Furthermore, we assume that \eqref{linear} has no bound states:
\eq
P_{sc}=1.\label{Pscasp}
\eeq
\begin{assumption}\label{dasp2}\eqref{Pscasp} holds.
\end{assumption}
According to \cite{Sof-W5}, under Assumptions \ref{dasp1} and \ref{dasp2}, the solution to \eqref{linear} satisfies Strichartz estimates as follows: 
\eq
\| U_{V}(t,0)\psi_0\|_{\s^q_t\s^r_x(\R^{5+1})}\lesssim \| \psi_0\|_{\s^2_x(\R^5)},\label{linear: Strichartz}
\eeq
for all $\psi_0\in \s^2_x(\R^5)$, where this holds for all real numbers $q$ and $r$, which satisfy the following conditions:
\begin{equation}
2 \leq q,r \leq \infty ; \quad \frac{2}{q} + \frac{5}{r} = \frac{5}{2}.\label{qrstr}
\end{equation}
Let us now consider a nonlinear problem. We take $V_{lin}(x,t)=V_0(x)+V_{s}(x,t)+V_c(x,t)$ and add a defocusing nonlinearity, $\lambda|\psi|$ ($\lambda>0$),  into the interaction. This leads to the following nonlinear system:
\eq
\begin{cases}
    i\partial_t\psi(x,t)=(-\Delta_x+V_{lin}(x,t)+\lambda|\psi|)\psi(x,t)\\
    \psi(x,0)=\psi_0\in \s^2_{rad}(\R^5)
\end{cases}, \quad (x,t)\in \R^{5+1}.\label{Nlinear}
\eeq
\begin{assumption}\label{dasp3}Assume that $V_{lin}(x,t)\in \s^\infty_tH^1_x(\R^{5+1})$ and is radial in $x$. Assume that $\psi_0$ leads to the following constraint:
\eq
\sup\limits_{t\in \R}\|\psi(t)\|_{H^1_x(\R^5)}<\infty.
\eeq
\end{assumption}
Under Assumptions \ref{dasp1}, \ref{dasp2} and \ref{dasp3}, we obtain 
$$
V_{lin}(x,t)+\lambda|\psi(t)|\in \s^\infty_tH^1_x(\R^{5+1}), 
$$
and as a consequence, according to Theorem 2.1 in \cite{SW20221}, we have
\eq
\psi_+:=s\text{-}\lim\limits_{t\to \infty} e^{itH_0}F_c(\frac{|x-2Pt|}{t^\alpha}\leq 1)\psi(t)\text{ exists in }H^1_x(\R^5).
\eeq
\begin{assumption}\label{dasp4}If Assumption \ref{dasp3} holds, then the solution to system \eqref{Nlinear}, $\psi(t)$, satisfies that 
\eq
\|\psi(t)-e^{-itH_0}\psi_+\|_{H^1_x(\R^5)}\to 0\label{smallness}
\eeq
as $t\to \infty$. 
\end{assumption}
The following proposition tells us that under Assumptions \ref{dasp1} to \ref{dasp4}, $\psi(t)$ satisfies Strichartz estimates:
\eq
\|\psi(t)\|_{\s^q_t\s^r_x(\R^{5+1})}\leq C\|\psi(0)\|_{\s^2_x(\R^5)}\label{Goal: Stri}
\eeq
for some constant $C=C(\sup\limits_{t\geq 0}\|\psi(t)\|_{H^1_x(\R^5)})>0$, where $q$ and $r$ satisfy \eqref{qrstr}:
\begin{proposition}Assuming that Assumptions \ref{dasp1} to \ref{dasp4} hold, \eqref{Goal: Stri} holds.
\end{proposition}
\begin{proof}
It suffices to build end-point Strichartz estimate:
\eq
\|\psi(t)\|_{\s^2_t\s^{10/3}_x(\R^{5+1})}\lesssim \|\psi(0)\|_{\s^2_x(\R^5)}.
\eeq
    Using Duhamel's formula with respect to $U_{V}(t,0)$, we obtain 
    \begin{align}
        \psi(t)=&U_{V}(t,0)\psi(0)+(-i)\lambda\int_0^tds U_{V}(t,s)|\psi(s)|e^{-isH_0}\psi_+\nonumber\\
        &+(-i)\lambda\int_0^tds U_{V}(t,s)|\psi(s)|\psi_r(s),
    \end{align}
where
\eq
\psi_r(t):=\psi(t)-e^{-itH_0}\psi_+.
\eeq
Set
\eq
C_V:=\|U_{V}(t,0)\|_{\s^2_x(\R^5)\to \s^2_t\s^{10/3}_x(\R^{5+1})}.
\eeq
$\psi_r(t)\in H^1_x(\R^5)$ implies $\psi_r(t)\in \s^{5/2}_x(\R^5)$ with
\eq
\|\psi_r(t)\|_{\s^{5/2}_x(\R^5)}\leq C\|\psi_r(t)\|_{H^1_x(\R^5)}
\eeq
for some constant $C>0$. Due to Assumptions \ref{dasp3} and \ref{dasp4}, there exists $T>0$ such that when $t\geq T$, the following inequality holds
\eq
\|\psi_r(t)\|_{\s^{5/2}_x(\R^5)}\leq C\|\psi_r(t)\|_{H^1_x(\R^5)}<\frac{1}{2C_V}\label{rsmall}.
\eeq
Using the inhomogeneous Strichartz estimate (see e.g., \cite{T2006})
\eq
\| \int_{s<t} dse^{-i(t-s)H_0}F(s)\|_{\s^q_t\s^r_x(\mathbb{R}^{n+1})}\lesssim_{n,q,r} \| F \|_{\s^{q'}_t\s^{r'}_x(\mathbb{R}^{n+1})},
\eeq
where $\frac{n}{r}+\frac{2}{q}=\frac{n}{2}, \text{ for }2\leq r,q\leq \infty,$ with $(q,r,n)\neq (2,\infty,2)$, and $(r,r')$ and $(q,q')$ are conjugate pairs, we have 
\begin{align}
    \|\psi(t)\|_{\s^2_t\s^{10/3}_x(\R^{5+1})}\leq& \| U_V(t,0)\psi(0)\|_{\s^2_t\s^{10/3}_x(\R^{5+1})}+C_V\| |\psi(t)|e^{-itH_0}\psi_+\|_{\s^2_t\s^{10/7}_x(\R^{5+1})}\nonumber\\
    &+C_V\| |\psi(t)|\psi_r(t)\|_{\s^2_t\s^{10/7}_x(\R^{5+1})}\nonumber \\
    \leq& \| U_V(t,0)\psi(0)\|_{\s^2_t\s^{10/3}_x(\R^{5+1})}+C_V\| |\psi(t)|e^{-itH_0}\psi_+\|_{\s^2_t\s^{10/7}_x(\R^{5+1})}\nonumber\\
    &+\frac{1}{2}\|\psi(t)\|_{\s^2_t\s^{10/3}_x(\R^{5+1})}+C_V\| \chi(t\in [0,T])|\psi(t)|\psi_r(t)\|_{\s^2_t\s^{10/7}_x(\R^{5+1})},\label{boot}
\end{align}
where we also use, due to \eqref{rsmall}, 
\begin{align}
    C_V\| \chi(t\geq T)|\psi(t)|\psi_r(t)\|_{\s^2_t\s^{10/7}_x(\R^{5+1})}\leq & C_V\|\psi_r(t)\|_{\s^\infty_t\s^{5/2}_x(\R^{5+1})}\|\psi(t)\|_{\s^2_t\s^{10/3}_x(\R^{5+1})}\nonumber\\
    \leq&\frac{1}{2}\|\psi(t)\|_{\s^2_t\s^{10/3}_x(\R^{5+1})}.
\end{align}
Let
\eq
E:=\sup\limits_{t\in \R}\|\psi(t)\|_{H^1_x(\R^5)}.
\eeq
Using \eqref{boot}, \eqref{linear: Strichartz} and H\"older's inequality, we obtain
\begin{align}
    \| \psi(t)\|_{\s^2_t\s^{10/3}_x(\R^{5+1})}\leq& 2\| U_V(t,0)\psi(0)\|_{\s^2_t\s^{10/3}_x(\R^{5+1})}+2C_V\| |\psi(t)|e^{-itH_0}\psi_+\|_{\s^2_t\s^{10/7}_x(\R^{5+1})}\nonumber\\
    &+2C_V\| \chi(t\in [0,T])|\psi(t)|\psi_r(t)\|_{\s^2_t\s^{10/7}_x(\R^{5+1})}\nonumber\\
    \lesssim& 2\| U_V(t,0)\psi(0)\|_{\s^2_t\s^{10/3}_x(\R^{5+1})}+2C_V\| |\psi(t)\|_{\s^\infty_t\s^{5/2}_x(\R^{5+1})}\|e^{-itH_0}\psi_+\|_{\s^2_t\s^{10/3}_x(\R^{5+1})}\nonumber\\
    &+2C_V\| \chi(t\in [0,T])|\psi(t)|\psi_r(t)\|_{\s^2_t\s^{10/7}_x(\R^{5+1})}\nonumber\\
    \lesssim_E& \|\psi(0)\|_{\s^2_x(\R^5)}+\sqrt{T}.
\end{align}
Here, we also use 
\eq
\| \psi_+\|_{\s^2_x(\R^5)}=\|\Omega_\alpha^*\psi(0)\|_{\s^2_x(\R^5)}\leq \|\psi(0)\|_{\s^2_x(\R^5)},
\eeq
\begin{align}
\|\psi_r(t)\|_{H^1_x(\R^5)}\leq& \|\psi(t)-e^{-itH_0}\Omega_\alpha^*\psi(0)\|_{H^1_x(\R^5)}\\
\lesssim & E,
\end{align}
and
\begin{align}
\||\psi(t)|\psi_r(t)\|_{\s^\infty_t\s^{10/7}_x(\R^{5+1})}\lesssim& \|\psi(t)\|_{\s^\infty_t\s^{20/7}_x(\R^{5+1})}\|\psi_r(t)\|_{\s^\infty_t\s^{20/7}_x(\R^{5+1})}\nonumber\\
\lesssim_E &1.
\end{align}
This completes the proof.
\end{proof}

\section{Epilogue}
In this section we discuss the issue of the precise meaning of Soliton Resolution.
Naively, it is the statement that all localized asymptotic states solve the time-independent NLS, giving a solution with a profile of a soliton.
However, this does not hold in general. It is known that asymptotic solutions of nonlinear equations that are not free waves, can be time dependent. Many such solutions have been found, and they are commonly referred to as {\bf coherent structures}, see e.g.\cites{Sof1, PW2022}. It includes breathers, kinks, monopoles, hedgehogs, black-holes, self-similar solutions, vortices,
peakons...
For the case of NLS type equations with potential terms, it is not clear what can come out. In \cite{Sof-W3} we showed that if the interaction term decays like $r^{-2}$ for large distance, we can construct global, stable, self-similar solutions that spread like $t^{\alpha},$ $\alpha<1/2,$ $n\geq 5.$
In one dimension, NLS systems of equations lead to breather solutions and consequently to time dependent potentials with time dependent bound states \cite{MSW-1}.
So the correct statement should be that soliton can be any coherent structure.\par

But there is more to the story: One in fact expects that many of those coherent structures to be unstable and non-generic.
In particular, we expect breathers for nonlinear dispersive and hyperbolic equations in the continuum, to be unstable and non-generic.
In fact we expect that excited state solitons are also unstable. See e.g. \cites{S-Wei5, S-Wei1}, \cite{Ts2002}.
 A major step in proving such a result, would be to show that time quasi-periodic solutions which are localized in space (and smooth) be unstable.
Now, we note that the results described in this paper, combined with the previous works of Tao \cite{T2004}, \cite{T2007}, \cite{T2008}, \cite{T2006}, \cite{T2014} and Liu-Soffer \cite{LS2020}, \cite{LS2021} imply that the localized solutions are in fact smooth to high order, depending only on the structure of the Interaction term.
 The proof of such instability of time quasi-periodic solutions should not be hard; here we describe the steps to proving this, based on known techniques.
The first step is to turn this question into a problem in spectral theory of a linear equation!\par

Suppose $\psi(x,t)$ is a localized smooth solution of the NLS type equation we discussed, which is moreover, time quasi-periodic.
This means that the following equation holds:
\begin{equation}
i\frac{\partial \psi}{\partial t}=-\Delta \psi- F(|\psi|)\psi.
\end{equation}
We now view $F(|\psi|)=V(x,t)=V(x,\omega_1 t,...\omega_N t)$ as a potential for a linear Schr\"odinger equation. We introduce a new Hilbert space where the quasi-periodicity is used:
$$
\mathcal{H}=L^2(\mathbb{R}^d\times [0,2\pi]\times...[0,2\pi]).
$$
On  this Hilbert space, we define the Hamiltonian
$$
L=-\Delta_x-\sum_{j=1}^{N}i\omega_j\frac{\partial }{\partial_{t_j}} +V(x,t_1,...t_N).
$$
Since the spectrum of $-i\frac{\partial }{\partial_{t_j}}$ consists of all the integers (times a constant),
the continuous spectrum of $L$ covers all of the real line. This is due to the fact the continuous spectrum of the Laplacian is the positive real line, and this will be unchanged by the perturbation, since it is relatively compact, due to the localization in space.

So, now we have a function $u$ that belongs to the Hilbert space $\mathcal{H}$ and is an eigenvector of $L.$
Hence, $L$ has an {\bf embedded eigenvalue} in its continuous spectrum.

The weakly localized part has some more important properties: It consists of  smooth functions of $x$ and therefore also of time.
Assuming, for a moment, that such solutions satisfy the Petite Conjecture (which will be discussed later), we can focus on the case where the time dependence is almost periodic.
The goal is to prove that such solutions are unstable and therefore non-generic.
That would leave us with the standard soliton as the only generic asymptotic state.\par

To proceed with this approach we first consider the case of a time periodic smooth and localized solution.

 The periodic solution corresponds to an $L^2$ eigenfunction of the extended operator $L$
$$
i\frac{\partial \psi}{\partial t}=-\Delta\psi +F(|\psi|)\psi,
$$
$$
-i\frac{\partial \psi}{\partial t}-\Delta\psi +F(|\psi|)\psi\equiv L\psi= E\psi.
$$
Here we assume that the interaction term is a function of $\psi$ only, but we could add a potential term.
We now observe that by basic spectral theory the continuous spectrum of $L$ on the Hilbert space of functions on $\mathbb{R}^d \times [0,T]$ covers the whole real line, with many overlapping branches.
In particular, it follows that the localized solution, if it exists, is an eigenfunction with an embedded eigenvalue.

Next, we appeal to the general theory of perturbations of embedded eigenvalues to conclude that such an eigenvalue is removed by a generic perturbation. Generic is in the sense that there is a number, called Fermi Golden Rule; when it is not zero, it implies the absence of the eigenvalue.\par
We will follow techniques from Linear Spectral Theory, in particular Resonance Theory.

Originally, this was proved using the assumption of Dilation analyticity of the potential. This is not suitable for our purposes. Another approach, more powerful and general, was introduced by Mourre \cite{cycon2009schrodinger}. This approach which basically reduces the problem to proving an appropriate Mourre estimate, also allows proving exponential decay of eigenfunctions in many cases. This was in fact used by Sigal \cite{ Sig} and Pyke-Sigal \cites{R-Sig, Pyke} to prove the absence of quasi-periodic solutions for both NLS and NLKG equations.
Some of the more general results are limited to the case of "small breathers", that is quasi-periodic solutions which are small amplitude bifurcations from an isolated eigenfunction of a linear Hamiltonian.

For our problem we can not use smallness assumption, and we propose to apply 
 the \emph{time dependent resonance theory}, originally developed by Soffer-Weinstein \cite{S-Wei3}, \cite{S-Wei4}, \cite{S-Wei2}.
In this approach only minimal decay of the potential is needed, and somewhat improved local decay estimate on the continuous spectral part of the Hamiltonian.
These conditions are further relaxed by Costin-Soffer \cite{CS2001}, to require that the resolvent of the Hamiltonian is a $C^{1+\eta}$ regular in the constant $E$ around the energy of the bound state. $\eta>0$ is sufficient.
Then the embedded eigenvalue is removed under any small perturbation with the Fermi Golden Rule satisfied. So, generically, for a given frequency of the periodic solution there is only one solution, or non. There is a complication in the use of resonance theory to the operator $L$ in the quasi-periodic case.
This is due to the fact that the spectral properties of $L$ are needed in an arbitrary small neighborhood of $0.$
Therefore, on the Hilbert subspaces where the spectrum of the operator $\sum_{j=1}^{N}i\omega_j\frac{\partial }{\partial_{t_j}}$ is close to zero, then we are in a threshold situation for the operator $L$ restricted to this subspace. 
The restrictions to such subspaces of $L$ looks like $-\Delta+V_0(x)+V_m(x,t)$ where $V_0(x)$ come from the time independent part of the potential, and $V_m$ come from very high oscillatory part of the potential. Due to the assumed smoothness such $V_m$ are very small.
So, we only need to have sufficient local decay for $-\Delta +V_0.$
In $5$ or more dimensions we have sufficient decay estimates in the generic case, which includes the case $V_0=0.$ \par
Next we consider the case of an almost periodic function.
In this case, the corresponding Hilbert space is non-separable, and one would need to extend spectral theory to such a space, which is not an easy proposition.
Instead, we propose to use the method of \emph{multiscale time averaging} \cite{FS}. In this approach, given a large time $T$ and small $\epsilon,$ one can average over a finite time interval the part of the potential with sufficiently small time frequencies, and get an averaged potential without low frequencies. The new Hamiltonian will approximate the full dynamics on an interval of order $T$ in time, with accuracy $\epsilon .$ Similar reduction can be used to remove all the high frequency oscillations.
Therefore since bound states are non-generic to the quasi-periodic approximation, the same is true for $L$, by choosing $T$ large enough.\par

Finally, we discuss the Petite Conjecture \cite{Sof1}.
This conjecture is the nonlinear analog of Ruelle's  geometric characterization of bound states.
One way to prove it is to rely on a Theorem of Wiener, that states that the Fourier transform $F(t)$ of a finite signed measure
satisfies 
$$
\lim_{T \to \infty} \frac{1}{T}\int_{0}^{T}|F(t)|^2 dt= \sum_j |\mu\{x_j\}|^2.
$$
This means that the discrete part of the spectrum is given by a sequence of (masses of) $\delta$ functions that converges in $L^2.$
Another useful way, also due to Wiener, is to prove that the autocorrelation function of $F(t)$ exists for all $t$ and is continuous in $t$.
Then such an autocorrelation function is positive definite, and it follows that it is the Fourier transform of a finite measure.\par
 
\appendix
\section{Appendix}
Here, we present a proof of the existence of free channel wave operators based on \cite{SW20221}.\par Here we consider a general class of nonlinear Schr\"odinger equations with $H_0:=-\Delta_x$ and a fixed parameter $a\in [0,1]$. The equation is of the form:
\eq
\begin{cases}
i\partial_t\psi-H_0\psi=\mathcal{N}(x,t,\psi)\psi \\
\psi(x,0)=\psi_0\in \mathcal{H}^a_x(\mathbb{R}^n)
\end{cases}, \quad (x,t)\in \mathbb{R}^n\times \mathbb{R}\label{A:SE}
\eeq
where $n\geq 1$ denotes the spatial dimension, and $\mathcal{N}(x,t,\psi)$ can take on any of the following forms: $V(x,t)$,  or $V(x,t)\psi+N(x,t,\psi)$. Note that $V(x,t)$ and $N(x,t,\psi)$ are not necessarily real. In \cite{SW20221}, we assume that the solution to system \eqref{A:SE} possesses a uniform $\mathcal{H}^a_x$ bound, which can be expressed as:
\eq
\sup\limits_{t\in \mathbb{R}}\|\psi(t)\|_{\mathcal{H}^a_x}\lesssim_{\|\psi_0\|_{\mathcal{H}^a_x}} 1\label{A:con: H}.
\eeq
$\mathcal{H}^a_x=\mathcal{H}^a_x(\mathbb{R}^n)$ denotes the Sobolev space.  The interaction terms $\mathcal{N}(x,t,|\psi|)$ can be one of the following:
\begin{enumerate}
\item \label{A:Lpresult}($\s^p$ potentials) $V(x,t)\in \s^\infty_t\mathcal{H}^a_x(\mathbb{R}^n\times\mathbb{R}), n\geq 3$. 

\item \label{A:Nonlinear}(nonlinear potentials) for some $a\in [0,1]$, we have $\mathcal{N}(x,t,\psi)=N(|\psi|)$ such that:
\eq
\|N(|\psi|)\|_{\mathcal{H}^a_x}\leq C(\|\psi\|_{\mathcal{H}^a_x}),
\eeq
and assuming that $\psi_0\in \mathcal{H}^a_x$ leads to:
\eq
\sup\limits_{t\in \mathbb{R}}\|\psi(t)\|_{\mathcal{H}^a_x}\lesssim_{\psi_0} 1.
\eeq
\end{enumerate}
Here, we denote $\langle \cdot\rangle: \mathbb{R}^n\to \mathbb{R}, x\mapsto \sqrt{|x|^2+1}$. Typical examples of nonlinear potentials are: 
\eq
\mathcal{N}(\psi)=P(|\psi|),\quad \text{when }|z|\leq 1, \quad |P(z)|\leq |z|^b,
\eeq
for $b>\max\{1,4/n\}$ with $P(z)$,
\begin{align}
&(n=3)\quad \mathcal{N}(\psi)=\pm [\frac{1}{|x|^{3/2-\delta}}*|\psi|^2](x), \quad \delta\in (0,\frac{3}{2}),
\end{align}
and 
\begin{align}
&(n=3)\quad \mathcal{N}(\psi)=\pm |\psi| \quad \text{Strauss Exponent}.
\end{align}
Note that the Strauss exponent corresponds to the case where $n=3$ and $\mathcal{N}(x,t,\psi)=\pm |\psi|$. We define $\s^p_{\delta,x}:={f(x): \langle x\rangle^\delta f(x)\in \s^p_x(\mathbb{R}^n)}$ for $1\leq p\leq\infty$. We introduce the smooth characteristic functions $\bar{\F}_c(\lambda)$ and $\F_j(\lambda)$ ($j=1,2$) of the interval $[1,+\infty)$, and define
\begin{equation}
\F_c(\lambda\leq a):=1-\bar{\F}_c(\lambda/a), \quad \F_j(\lambda>a):=\F_j(\lambda/a), \quad j=1,2,
\end{equation}
as well as
\begin{equation}
\bar{\F}_c(\lambda\leq a):=\bar{\F}_c(\lambda/a), \quad \bar{\F}_j(\lambda\leq a):=1-\F_j(\lambda/a),
\end{equation}
where for $j=1,2,$ $F_j(k)=\bar{F}_c(k)=1$ when $k\geq 1$ and $F_j(k)=\bar{F}_c(k)=0$ when $k\leq 1/2$.
We now present our main results:
\begin{theorem}\label{A:thm1}Let $\psi(t)$ be a global solution of equation \eqref{A:SE} with an interaction satisfying \eqref{A:Lpresult}, i.e., $V(x,t)=N(x,t,|\psi|)\in \s^\infty_t\mathcal{H}^a_x(\mathbb{R}^n\times\mathbb{R})$ for $n\geq 3$ and some $a\in [0,1]$. For $\alpha\in(0,1-2/n)$, the channel wave operator acting on $\psi(0)$, given by
\eq
\Omega_\alpha^*\psi(0)=s\text{-}\lim\limits_{t\to \infty} e^{itH_0}\F_c(\frac{|x-2tP|}{t^\alpha}\leq 1)\psi(t), \label{A:wave-4}
\eeq
exists in $\mathcal{H}^a_x(\R^n)$, and 
\eq
w\text{-}\lim\limits_{t\to \infty} e^{itH_0}\bar{\F}_c(\frac{|x-2tP|}{t^\alpha}> 1)\psi(t)=0,\quad \text{ in }\mathcal{H}^a_x(\R^n).\label{A:weakthm1}
\eeq
Therefore, we have
\eq
\|\psi(t)-e^{-itH_0}\Omega_{\alpha}^*\psi(0)-\psi_{w}(t)\|_{\mathcal{H}^a_x}=0
\eeq
for all $t\geq 0$, where 
\eq
\psi_{w}(t):=\bar{\F}_c(\frac{|x-2tP|}{t^\alpha}>1)\psi(t).
\eeq
\end{theorem}
\begin{remark}
The condition on the interaction here is not sharp. The same proof applies to an abstract version \cite{SW20221}, which is more general.
\end{remark}
\begin{remark}To control the non-free part, i.e., the weakly localized part, even in three or higher dimensions, we require that the interaction term is localized in the spatial variable $x$. Note that the interaction can also be nonlinear, including monomial terms, thanks to the use of radial Sobolev embedding theorems for $H^1$ functions in three or higher dimensions.
\end{remark}
\section{Propagation Estimates And Proof of Theorem \ref{A:thm1}}
\subsection{Propagation Estimate}
Given an operator family $B(t)$, we denote
\eq
\langle B\rangle_t:=(\psi(t), B(t)\psi(t))_{\s^2_x(\mathbb{R}^n)}=\int_{\mathbb{R}^n} \psi(t)^*B(t)\psi(t)d^nx.
\eeq 
Suppose a family of self-adjoint operators $B(t)$ satisfy the following estimate:
\begin{align}
&\partial_t\langle B\rangle_t=(\psi(t), C^*C\psi(t))_{\s^2_x(\mathbb{R}^n)}+g(t),\nonumber\\
&g(t)\in L^1(dt), \quad C^*C\geq 0.
\end{align}
In this case, we refer to the family $B(t)$ as a \textbf{Propagation Observable} (PROB) [\cite{HSS1999}, \cite{SS1988}, \cite{SS1987}].

Upon integration over time, we obtain the bound:
\begin{subequations}\label{A:CC}
\begin{align}
\int_{t_0}^T\|C(t)\phi(t) \|_{\s^2_x(\mathbb{R}^n)}^2dt=&( \psi(T), B(T)\psi(T))_{\s^2_x(\mathbb{R}^n)}-( \psi(t_0), B(t_0)\psi(t_0))_{\s^2_x(\mathbb{R}^n)}-\int_{t_0}^Tg(s) ds\nonumber\\
\leq& \sup\limits_{t\in [t_0,T]} \left|( \psi(t), B(t)\psi(t))_{\s^2_x(\mathbb{R}^n)}\right|+C_g, 
\end{align}
where
\eq
C_g:=\|g(t)\|_{\s^1_t(\mathbb{R})}.
\eeq
\end{subequations}
\subsection{Relative Propagation Estimate}
In this paper, we present a modified version of the {\bf PRES}, which we refer to as the Relative {\bf PRES} (RPRES). Given an operator $\tilde{B}$, we define its time-dependent expectation value as
\eq
\langle \tilde{B}: \phi(t)\rangle_t:=(\phi(t), \tilde{B}(t)\phi(t)  )_{\s^2_x(\mathbb{R}^n)}=\int_{\mathbb{R}^n}\phi(t)^*\tilde{B}(t)\phi(t)d^nx,
\eeq
where $\phi(t)$ does not need to be the solution to \eqref{A:SE}, but it satisfies
\eq
\sup\limits_{t\geq 0}\langle \tilde{B}:\phi(t)\rangle_t<\infty.\label{phiH}
\eeq
Suppose \eqref{phiH} is satisfied, and $\partial_t\langle \tilde{B}:\phi(t)\rangle_t$ satisfies the following estimate:
\begin{align}
&\partial_t\langle \tilde{B} : \phi(t)\rangle_t=\pm\langle \phi(t), C^*C\phi(t)\rangle+g(t)\nonumber\\
&g(t)\in L^1(dt), \quad C^*C\geq 0 .
\end{align}
We then refer to the family $\tilde{B}(t)$ as a {\bf Relative Propagation Observable}(RPROB) with respect to $\phi(t)$.

Upon integration over time, we obtain the bound:
\eq
\int_{t_0}^T\|C(t)\phi(t) \|_{\s^2_x(\mathbb{R}^n)}^2dt\leq \sup\limits_{t\in [t_0,T]} \left|( \phi(t), \tilde{B}(t)\phi(t))_{\s^2_x(\mathbb{R}^n)}\right|+C_g, \quad C_g:=\|g(t)\|_{\s^1_t(\mathbb{R})}.\label{A:CC0}
\eeq
In this appendix, we choose $\phi(t)=e^{itH_0}\psi(t)$. The operator $C^*C(t)$ is a multiplication operator $\partial_t[F_c](\frac{|x|}{t^\alpha}\leq 1)\in \s^1_t[1,\infty)$ when $n\geq 3$, or a multiplication operator in frequency space $\partial_t[F_1](t^b|q|\geq 1)\in \s^1_t[1,\infty)$ (where $q$ denotes the frequency variable) when $n\geq 1$. Then by using H\"older's inequality, \eqref{A:CC} implies that for $T\geq t_0\geq 1$,
\eq
\| \int_{t_0}^T dt C^*(t)C(t)\psi(t)\|_{\s^2_x}\leq \left(\int_{t_0}^T dt |C(t)|^2 \right)^{1/2}\left(\int_{t_0}^T\|C(t)\psi(t) \|_{\s^2_x}^2dt\right)^{1/2} \to 0
\eeq
as $t_0\to \infty$.
We call this estimate {\bf Relative Propagation Estimate}(RPRES).

We conclude this subsection by presenting an abstract formulation of the main theorem regarding the Free Channel Wave Operator.
\begin{theorem}\label{abstract} 
Let $H_0=\omega(p)$ be the generator of the free flow $U_0(t,0)=e^{-iH_0 t}$ acting on a Hilbert space $\mathcal{H}=L_x^2(\R^n).$
Let $U(t,0)$ be the solution operator of the Schr\"odinger type equation
$$
i\frac{\partial \psi(t)}{\partial t}=(H_0+V(x,t))\psi(t),
$$
where $V(x,t)=N(x,t,|\psi|).$ Assume that for initial data $\psi_0$ the solution of the above (nonlinear) equation is global, uniformly bounded in the Sobolev space $H_x^1(\R^n).$
Assume, moreover that the group $U_0$ is bounded from $H_x^{s,p}(\R^n)$ into $L_x^{p'}(\R^n)$ with a bound that decays faster than $1/t^{1+\epsilon}$ for some $ \epsilon >0.$
Here $p<2, s\geq 0.$

Then, if the above conditions are satisfied with $s=0,$   the following strong limit, defining the Free Channel Wave Operator exists:
\begin{align}[Free~ Channel]~~\quad\label{A:Free}
&\Omega_F(\psi_0)\equiv s\text{-}\lim_{ t\to \infty} U_0(0,t)J_c\psi(t)\\
&\psi(t)=U(t,0)\psi_0, \\
&J_c= U_0(t,0)J(|x|/t^{\alpha} \leq 1)U_0(0,t),\\
&0\leq \alpha\leq a(n,p)<1,
\end{align}
provided the interaction term satisfies the following bound uniformly in time:
\begin{equation}\label{A:interaction}
\| V(x,t)\psi(t)\|_{L^p_x(\R^n)} \lesssim 1.
\end{equation}
If the above assumptions hold with $s>0,$ we have that the same limit holds, provided we replace $J_c$ above by
$$
J_c= U_0(t,0)F_2(|p|/t^{\beta} \leq 1)J(|x|/t^{\alpha} \leq 1)F_2(|p|/t^{\beta} \leq 1)U_0(0,t)
$$
with $\beta$ sufficiently small compared with $s$ and $\epsilon$ mentioned above.

 Assume moreover that the interaction term is localized in space, uniformly in time, with decay rate $\delta>1.$
 Then, by modifying the channel wave operator using the change of the original $J_c$ to  $$U_0(t,0)\Sigma F_1(|P-\tau|>t^{-1/2+\beta})J(|x|/t^{\alpha}\leq 1)\Sigma F_1(|P-\tau|>t^{-1/2+\beta}) U_0(0,t),$$  the channel wave operators exist, provided local decay estimates hold away from the thresholds of the operator $\omega(P).$
Here, $ \Sigma F_1(|P-\tau|>t^{-1/2+\beta})$ stands for the sum of cutoff functions over all thresholds $\tau=\{\tau'| \nabla_k \omega(k)=0 \, \text{for}\,        k=\tau'\}$ of $\omega.$

\end{theorem}
\begin{remark}
 The limit is zero on the support of $1-F_2$ due to the assumption that the solution is uniformly bounded in $H_x^1(\mathbb{R}^n).$
\end{remark}
\subsection{Time translated Potential}
Given a potential $V$, the {\bf time-translated ($tT$) potential} is defined as the translation of $V$ under the flow of the free Hamiltonian $H_0:=-\Delta_x$, i.e.,
\begin{equation}
\sK_t(V) := e^{itH_0}Ve^{-itH_0},
\end{equation}
as introduced in \cite{SW2020}. The $tT$ potential has the following representation formulas:
\begin{align}
\sK_t(V) &= \frac{1}{(2\pi)^{n/2}}\int d^n\xi, \hat{V}(\xi,t) e^{i(x+2tP)\cdot \xi}, \\
\sK_t(V) &= V(x+2tP, t). \label{A:sKre1}
\end{align}
Here, $P:=-i\nabla_x$, and $\hat{V}(\xi,t)$ denotes the Fourier transform of $V(x,t)$ with respect to the $x$ variables:
\begin{equation}
\hat{V}(\xi,t) := \frac{1}{(2\pi)^{n/2}}\int d^nx, e^{-ix\cdot \xi} V(x,t).
\end{equation}
We will also use $c_n$ to represent $1/(2\pi)^{n/2}$.
\subsection{Estimates of the interaction terms and proof of Theorem \ref{A:thm1}}When the space dimension $n\geq 3$, we define
\eq
\Omega_\alpha^*(t)\psi(0):=\langle P\rangle^{-a}e^{itH_0}\mathcal{F}_c(\frac{|x-2tP|}{t^\alpha}\leq 1)\langle P\rangle^a\psi(t)\quad \text{ in }\mathcal{H}_x^a(\R^n)
\eeq
for $\psi(0)\in \mathcal{H}^a_x(\R^n)$. Here, we have $\Omega_\alpha^*(t)\psi(0)\in \mathcal{H}_x^a$ since $\psi(t)\in \mathcal{H}_x^a$. We use $\langle P\rangle^{-a} \mathcal{F}_c\langle P\rangle^a$ instead of $\mathcal{F}_c$ for convenience, because
\eq
\|\langle P\rangle^{-a}e^{itH_0}\mathcal{F}c(\frac{|x-2tP|}{t^\alpha}\leq 1)\langle P\rangle^a\psi(t) \|_{\mathcal{H}_x^a}=\|e^{itH_0}\mathcal{F}c(\frac{|x-2tP|}{t^\alpha}\leq 1)\langle P\rangle^a\psi(t) \|_{L_x^2(\R^n)}
\eeq
and by using \eqref{A:sKre1},
\begin{align}
&s\text{-}\lim\limits_{t\to \infty} \Omega_\alpha^*(t)\psi(0)-e^{itH_0}\F_c(\frac{|x-2tP|}{t^\alpha}\leq 1)\psi(0)\nonumber\\
=&s\text{-}\lim\limits_{t\to \infty} (\langle P\rangle^{-a}\F_c(\frac{|x|}{t^\alpha}\leq 1)\langle P\rangle^a-\F_c(\frac{|x|}{t^\alpha}\leq 1))e^{itH_0}\psi(0)\nonumber\\
=&s\text{-}\lim\limits_{t\to \infty} \left[\langle P\rangle^{-a},\F_c(\frac{|x|}{t^\alpha}\leq 1)\right]\langle P\rangle^ae^{itH_0}\psi(0)\nonumber\\
=&0.\label{A:a<P>}
\end{align}
Using \eqref{A:sKre1} to rewrite $\Omega_\alpha^*(t)\psi(0)$ and applying {\bf Cook's method} to expand it, we get
\eq
\Omega_\alpha^*(t)\psi(0)=\langle P\rangle^{-a}\mathcal{F}_c(\frac{|x|}{t^\alpha}\leq 1)e^{itH_0}\langle P\rangle^a\psi(t),
\eeq
\begin{align}
\Omega_\alpha^*(t)\psi(0)=&\Omega_\alpha^*(1)\psi(0)+\int_1^tds  \langle P\rangle^{-a}  \partial_s[\mathcal{F}_c(\frac{|x|}{s^\alpha}\leq 1)]e^{isH_0}\langle P\rangle^a\psi(s)\nonumber\\
&+(-i)\int_1^tds \langle P\rangle^{-a} \mathcal{F}_c(\frac{|x|}{s^\alpha}\leq 1)e^{isH_0}\langle P\rangle^aV(x,s)\psi(s)\nonumber\\
=:&\Omega_\alpha^*(1)\psi(0)+\psi_p(t)+\int_1^tds\psi_{in}(s).
\end{align}
Due to \eqref{A:con: H}, 
\eq
\| \Omega_\alpha^*(1)\psi(0)\|_{\mathcal{H}_x^a}\lesssim \sup\limits_{t\in\R} \|\psi(t)\|_{\mathcal{H}^a_x}.
\eeq
We refer to $\psi_{in}(t)$ as the {\bf interaction term}. If $\|\psi_{in}(t)\|_{\mathcal{H}^a_x}\in L^1_t[1,\infty)$, then 
\eq
\int_1^\infty \psi_{in}(s)ds=\lim\limits_{t\to \infty}\int_1^t \psi_{in}(s)ds \text{ exists in }\mathcal{H}^a_x.\label{A:new psiin}
\eeq
Furthermore, by taking 
$$
\begin{cases}
B(t)=\F_c(\frac{|x|}{t^\alpha}\leq 1)\\
\phi(t)=e^{itH_0}\langle P\rangle^a\psi(t)
\end{cases}
$$
and using {\bf RPRES} with respect to $\phi(t)$, we obtain
\eq
\lim\limits_{t\to \infty}\psi_p(t)\text{ exists in }\mathcal{H}^a_x(\R^n):\label{A:new psip}
\eeq
\begin{subequations}\label{A:Cook}
   \eq
\langle B(t): \phi(t)\rangle_t\leq \left(\sup\limits_{t\in \R}\|\psi(t)\|_{\mathcal{H}_x^a}\right)^2.
\eeq
$\partial_t[\langle B(t):\phi(t)\rangle_t]$ reads
\eq
\partial_t[\langle B(t):\phi(t)\rangle_t]=c_p(t)+g(t),
\eeq
with 
$$
\partial_t[\F_c(\frac{|x|}{t^\alpha}\leq 1)]\geq 0,\quad \forall t\geq 1,
$$
$$
c_p(t):=(e^{itH_0}\langle P\rangle^a\psi(t), \partial_t[\F_c(\frac{|x|}{t^\alpha}\leq 1)]e^{itH_0}\langle P\rangle^a\psi(t))_{L^2_x(\R^n)}\geq 0,\quad \forall t\geq 1,
$$
and
\begin{align}
g(t):=&(-i)(e^{itH_0}\langle P\rangle^a\psi(t), \F_c(\frac{|x|}{t^\alpha}\leq 1)e^{itH_0}\langle P\rangle^a V(x,t)\psi(t))_{L^2_x(\R^n)}\nonumber\\
&+i(e^{itH_0}\langle P\rangle^aV(x,t)\psi(t), \F_c(\frac{|x|}{t^\alpha}\leq 1)e^{itH_0}\langle P\rangle^a\psi(t))_{L^2_x(\R^n)}.
\end{align}
Using $\|\psi_{in}(t)\|_{\mathcal{H}_x^a}\in L^1_t[1,\infty)$ and H\"older's inequality, we have
\eq
\|g(t)\|_{L^1_t[1,\infty)}\leq 2\left(\sup\limits_{t\in \R}\|\psi(t)\|_{\mathcal{H}_x^a}\right)\times \int_1^\infty \| \psi_{in}(t)\|_{\mathcal{H}_x^a}dt<\infty.
\eeq
Then for all $T\geq 1$,
\begin{align}
\int_1^T c_p(t)dt\leq &\langle B(t):\phi(t)\rangle_t\vert_{t=T}-\langle B(t):\phi(t)\rangle_t\vert_{t=1}+\|g(t)\|_{L^1_t[1,\infty)}\nonumber\\\leq & \left(\sup\limits_{t\in \R}\|\psi(t)\|_{\mathcal{H}_x^a}\right)^2+\|g(t)\|_{L^1_t[1,\infty)}.
\end{align}
Hence, 
\eq
\int_1^\infty c_p(t)dt=\lim\limits_{T\to \infty}\int_1^T c_p(t)dt\quad \text{ exists in }\R
\eeq
and for all $t_1>t_2\geq T\geq 1$, using H\"older's inequality in $s$ variable, 
\begin{align}
&\|\psi_p(t_1)-\psi_p(t_2)\|_{\mathcal{H}_x^a}=\|\int_{t_2}^{t_1}\partial_s[\F_c(\frac{|x|}{s^\alpha}\leq 1)]e^{isH_0}\langle P\rangle^a \psi(s) ds\|_{L^2_x(\R^n)}\nonumber\\
\leq & \|\left(\int_{t_2}^{t_1}\partial_s[\F_c(\frac{|x|}{s^\alpha}\leq 1)]ds\right)^{1/2} \left(\int_{t_2}^{t_1} \partial_s[\F_c(\frac{|x|}{s^\alpha}\leq 1)]|e^{isH_0}\langle P\rangle^a \psi(s)|^2ds\right)^{1/2} \|_{L^2_x(\R^n)}\nonumber\\
\leq &\left( \int_{T}^{\infty}c_p(s)ds\right)^{1/2}\to 0
\end{align}
as $T\to \infty$. So $\{\psi_p(t)\}_{t\geq 1}$ is Cauchy in $\mathcal{H}^a_x(\R^n)$ and \eqref{A:new psip} is true.
\end{subequations}
Therefore, we establish the existence of a free-channel wave operator on $\mathcal{H}^a_x$ due to \eqref{A:new psiin} and \eqref{A:new psip}, that is, 
\eq
s-\lim\limits_{t\to \infty} \Omega_\alpha^*(t)\psi(0)=\Omega_\alpha^*(1)\psi(0)+\psi_p(\infty)+\int_1^{\infty}ds\psi_{in}(s)\quad\text{ exists in }\mathcal{H}_x^a.
\eeq

For the interaction term, we also used the $L^p$ decay estimates of the free flow:
\eq
\|e^{-itH_0}f(x)\|_{\s^p_x(\mathbb{R}^n)}\lesssim_n \frac{1}{|t|^{\frac{n}{2}(\frac{1}{2}-\frac{1}{p})}}\|f(x)\|_{\s^{p'}_x(\mathbb{R}^n)},\quad f\in \s^{p'}_x(\mathbb{R}^n),
\eeq
where
\eq
\frac{1}{p}+\frac{1}{p'}=1, \quad 2\leq p\leq\infty.
\eeq

Throughout the paper, we need the following estimates for the interaction terms. 
\begin{definition}[Potential $V(x,t)$]Given $\psi(0)$ such that the solution $\psi(t)$ is bounded in $\mathcal{H}^a_x$ uniformly in time, we define $V(x,t)$ as $N(x,t,|\psi(t)|)$.
\end{definition}
We will also use $V(x,t)$ as a potential for a \emph{linear problem}. In particular, we remind the reader that if $V$ is assumed to be localized in $x$, this is achieved in the nonlinear case by assuming radial symmetry and applying the radial Sobolev Embedding Theorems (in three or more dimensions).
\begin{proposition}\label{A:can2}For $V(x,t)\in \s^\infty_t\mathcal{H}^a_{x}(\mathbb{R}^n\times \mathbb{R})$, $a\in[0,1]$, $\alpha\in(0,1-2/n)$, we have that for $t\geq 1, n\geq 3$,
\begin{multline}
\|\langle P\rangle^{-a}\F_c(\frac{|x|}{t^\alpha}\leq 1)e^{itH_0}\langle P\rangle^{a}V(x,t)\psi(t) \|_{\mathcal{H}_x^a}\lesssim_{n} \\
\frac{1}{t^{1+\beta}}\|V(x,t)\|_{\s^\infty_t\mathcal{H}^a_{x}(\mathbb{R}^n\times \mathbb{R})}\times\sup\limits_{t\in \R}\|\psi(t)\|_{\mathcal{H}^a_x(\mathbb{R}^n)}\in L^1_t[1,\infty),\label{A:can2eq}
\end{multline}
for some $\beta$ satisfying
\eq
\beta:=\frac{n(1-\alpha)}{2}-1>0.
\eeq
\end{proposition}
\begin{proof}
By using $\s^p$ decay of the free flow, we have
\begin{align}
\|\langle P\rangle^{-a}\F_c(\frac{|x|}{t^\alpha}\leq 1)e^{itH_0}&\langle P\rangle^{a}V(x,t)\psi(t) \|_{\mathcal{H}_x^a}=\|\F_c(\frac{|x|}{t^\alpha}\leq 1)e^{itH_0}\langle P\rangle^{a}V(x,t)\psi(t) \|_{L_x^2(\R^n)}\nonumber\\
\lesssim_n& \|\F_c(\frac{|x|}{t^\alpha}\leq 1) \|_{L^2_x(\R^n)}\times \frac{1}{t^{n/2}}\| V(x,t)\psi(t)\|_{W^{a,1}_x(\R^n)}\nonumber\\
\lesssim_n &\frac{1}{t^{n(1-\alpha)/2}}\|V(x,t)\|_{\s^\infty_t\mathcal{H}^a_x}\times\sup\limits_{t\in \R}\|\psi(t)\|_{\mathcal{H}^a_x(\mathbb{R}^n)}\in L^1_t[1,\infty)
\end{align}
provided that $\alpha\in (0,1-2/n)$.

\end{proof}
\begin{remark} Based on the proof of Proposition \ref{A:can2}, $L^\infty$ decay estimates of the free flow are not necessary in $n\geq 3$ dimensions. For example, $\s^{6+\epsilon}$ decay will be sufficient in $n=3$ dimensions.
\end{remark}
\begin{proof}[Proof of Theorem \ref{A:thm1}]
The existence of $\Omega_\alpha^*\psi(0)$ in $\mathcal{H}^a_x$ follows from \eqref{A:a<P>}, Proposition \ref{A:can2} and \eqref{A:new psip}, which is proved by using \eqref{A:Cook}.  \eqref{A:weakthm1} is based on the following argument: For all $\phi(x)\in \mathcal{H}^{-a}_x(\R^n)$, we can utilize \eqref{A:sKre1} and \eqref{A:con: H} to obtain
\begin{align}
    (\phi(x), e^{itH_0}\bar{\F}_c(\frac{|x-2tP|}{t^\alpha}>1)\psi(t))_{\s^2_x(\R^n)}=&(\langle P\rangle^{-a}\bar{\F}_c(\frac{|x|}{t^\alpha}>1)\phi(x), \langle P\rangle^ae^{itH_0}\psi(t))_{\s^2_x(\R^n)}\nonumber\\
    \to &0,
\end{align}
since
\eq
\| \langle P\rangle^{-a}\bar{\F}_c(\frac{|x|}{t^\alpha}>1)\phi(x)\|_{\mathcal{H}^{-a}_x(\R^n)}\to 0
\eeq
as $t\to \infty$. We have thus completed the proof.
\end{proof}

\noindent{\textbf{Acknowledgements}: }A. S. was partially supported by Simons Foundation Grant number 851844 and NSF grants DMS-2205931. X. W. was partially supported by NSF Grant DMS-1802170,  NSF Grant DMS-21-06255 and NSF Grant DMS-2205931 .

\bibliographystyle{amsplain}
\bibliography{bib}

\vspace{2mm}

\noindent
\textsc{Soffer: Department of Mathematics, Rutgers University, Piscataway, NJ 08854, U.S.A.}\\
{\em email: }\textsf{\bf soffer@math.rutgers.edu}

\medskip\noindent
\textsc{Wu: Department of Mathematics, Rutgers University, Piscataway, NJ 08854, U.S.A.}\\
{\em email: }\textsf{\bf xw292@math.rutgers.edu}


\end{document}